\documentclass[12pt]{amsart}
\usepackage{a4wide}
\usepackage{amssymb}
\newenvironment{pf}{\medskip\noindent{\it Proof.}\enspace}%
  {\hfill$\square$\newline\smallskip}
\newtheorem{lem}{Lemma}
\newtheorem{theo}[lem]{Theorem}
\newtheorem{cor}[lem]{Corollary}

\def\q{\hskip0.17cm}
\def\,{\hskip0.12cm}

\def\frak{\mathfrak}
\def\Bbb{\mathbb}
\def\cal{\mathcal}

\def\m{\begin{pmatrix}}
\def\em{\end{pmatrix}}
\def\sm{\left(\smallmatrix}
\def\esm{\endsmallmatrix\right)}

\begin{document}
\voffset=24pt
\title{\bf Super-replicable functions ${\cal N}(j_{1,N})$
           and periodically vanishing property}
\author{Chang Heon Kim \and Ja Kyung Koo}

\address{Chang Heon Kim, Department of mathematics, Seoul Women's university,
126 Kongnung 2-dong, Nowon-gu, Seoul, 139-774 Korea}
\email{chkim@swu.ac.kr}

\address{Korea Advanced Institute of Science and
Technology, Department of Mathematics, Taejon, 305-701 Korea}
\email{jkkoo@math.kaist.ac.kr}

 \begin{abstract}
 We find the super-replication formulae which would be a
 generalization of replication formulae. And we apply the formulae
 to derive periodically vanishing property in the Fourier
 coefficients of the Hauptmodul ${\cal N}(j_{1,12})$ as a
 super-replicable function.
 \end{abstract}

 \maketitle
 \renewcommand{\thefootnote}%
             {}
 \footnotetext{Supported by KOSEF Research Grant 98-0701-01-01-3
 \par AMS
  Classification : 11F03,
       11F22}
 \baselineskip=24pt

\section{Introduction}
\par Let ${\frak H}$ be the complex upper half plane and let $\Gamma_1(N)$
be a congruence subgroup of $SL_2({\Bbb Z})$ whose elements are
congruent to $\sm 1&*\\0&1 \esm \mod N \, (N=1,2,\cdots)$. Since
the group $\Gamma_1(N)$ acts on ${\frak H}$ by linear fractional
transformations, we get the modular curve
$X_1(N)=\Gamma_1(N)\backslash{\frak H}^*$, as a projective closure
of the smooth affine curve $\Gamma_1(N)\backslash{\frak H}$, with
genus $g_{1,N}$.
 Here, ${\frak H}^*$ denotes the union of ${\frak H}$ and $\Bbb
 P^1({\Bbb Q})$.
 \par Ishida and Ishii showed in \cite{Ishida-Ishii}
 that for $N\ge 7$, the function field $K(X_1(N))$ is generated
 over ${\Bbb C}$ by the modular functions $X_2(z,N)^{\epsilon_N
 \cdot N}$ and $X_3(z,N)^N$, where
 $X_r(z,N)=e^{2\pi i \frac{(r-1)(N-1)}{4N}} \prod_{s=0}^{N-1}
 \frac{K_{r,s}(z)}{K_{1,s}(z)}$ and $\epsilon_N$ is $1$ or $2$
 according as $N$ is odd or even. Here, $K_{r,s}(z)$ is a Klein
 form of level $N$ for integers $r$ and $s$ not both congruent to
 $0 \mod N$.
 On the other hand, since the genus $g_{1,N}=0$ only for the
eleven cases $1\le N \le 10$ and $N=12$ (\cite{Miyake},
\cite{Chang2}), the function field $K(X_1(N))$ in this case is a
rational function field ${\Bbb C}(j_{1,N})$ for some modular
function $j_{1,N}$ (Table 3, Appendix).
 \par The element $\sm 1&1 \\ 0&1 \esm$ of $\Gamma_1(N)$ takes
 $z$ to $z+1$, and in particular a modular function $f$ in
 $K(X_1(N))$ is periodic. Thus it can be written as a Laurent series in
 $q=e^{2\pi i z}$ ($z\in {\frak H}$), which is called a $q$-{\it series} \, (or
 $q$-{\it expansion}) \, of $f$. We call $f$ {\it normalized} \, if
 its $q$-series starts with $q^{-1}+0+a_1 q+a_2 q^2+\cdots$. By a
 {\it Hauptmodul} $t$ we mean the normalized generator of
 a genus zero function field $K(X_1(N))$ and
 we write $t=q^{-1}+0+\sum_{k\ge 1} H_k q^k $ for its $q$-series.
\par For a Fuchsian group $\Gamma$, let $\overline{\Gamma}$ denote
the inhomogeneous group of $\Gamma$ ($=\Gamma/\pm I$). Let
$\Gamma_0(N)$ be the Hecke subgroup given by $\{ \sm a&b\\c&d \esm
\in SL_2({\Bbb Z}) \mid c\equiv 0 \mod N \}$. Also, let $t={\cal
N}
 (j_{1,N})$ be the
Hauptmodul of $\Gamma_1(N)$ and $X_n(t)$ be a unique polynomial in
$t$ of degree $n$ such that $X_n(t)-\frac 1n q^{-n}$ belongs to
the maximal ideal of the local ring ${\Bbb C}[[q]]$. Polynomials
with this property are known as the Faber polynomials
(\cite{Duren}, Chapter 4). Write $X_n(t)=\frac 1n q^{-n} +
\underset{m\ge 1}{\sum} H_{m,n} q^m$.
 \par When
 $\overline{\Gamma}_1(N)=\overline{\Gamma}_0(N)$, ${\cal N}
 (j_{1,N})$ becomes a replicable function, that is, it satisfies
 the following replication formulae
 \begin{equation*}
  H_{a,b}=H_{c,d} \text{ \q \q
 whenever $ab=cd$ and
 $(a,b)=(c,d)$} \tag{$*$}
 \end{equation*} (\cite{ACMS}, \cite{Cummins-Norton}, \cite{Koike}).
  Given a replicable function $f$
  the {\it $n$-plicate} of $f$ is defined iteratively by
 $$ f^{(n)}(nz)= -\sideset{}{'}\sum_{\smallmatrix ad=n \\ 0\le b < d
 \endsmallmatrix}
   \, f^{(a)}\left(\frac{az+b}{d}\right)+nX_n(f) $$
 where the primed sum means that the term with $a=n$ is omitted (\cite{Cummins-Norton}).
  We call $f$ {\it completely replicable} if $f$ is a replicable
 function with rational integer coefficients and has only a finite
 number of distinct replicates, which are themselves replicable functions.
 According to \cite{ACMS} there are, excluding the trivial cases $q^{-1}+aq$,
 326 completely replicable functions of which 171 are monstrous
 functions, i.e.,  modular
 functions whose $q$-series coincide with the Thompson series
 $T_g(q)=\sum_{n\in {\Bbb Z}} {\rm Tr}(g| V_n)q^n$ for some element
 $g$ of the monster simple group $M$ whose order is approximately
 $8\cdot 10^{53}$. Here we observe that $V=\underset{n\in {\Bbb Z}}{\bigoplus} V_n$
 is the infinite dimensional graded representation of $M$
 constructed by Frenkel {\it et al.} (\cite{Frenkel1},
 \cite{Frenkel2}).
 Furthermore, in \cite{Cummins-Norton} Cummins and Norton showed that
 if $f$ is replicable, it can be
 determined only by the 12 coefficients of its first 23 ones.
 \par If $\overline{\Gamma}_1(N)\neq
 \overline{\Gamma}_0(N)$, unlike those replicable functions mentioned above,
 we show in \S 3 that the Fourier coefficients of $X_n(t)$ with $t={\cal N}(j_{1,N})$
 ($N \neq 7, 9$)
  satisfy a twisted
 formula (\ref{super}) by a character $\psi$ (see Corollary
 \ref{Repli2}). Here we note that when we work with the
 Thompson series, it is reduced to replication formulae in $(*)$ by
 viewing $\psi$ as the trivial character. Thus in this sense it gives
 a more general class of modular functions, which we
 propose to call ${\cal N}(j_{1,N})$ a {\it super-replicable}
 function.
 \par There would be certain similarity between some of
 replicable functions and super-replicable ones as follows.
 We derived in \cite{Kim-KooA}
 the following self-recursion formulas for the Fourier
 coefficients of ${\cal N}(j_{1,N})$ without the aid
 of its 2-plicate when
 $N=2,6,8,10,12$:  for $k\ge 1$,
\begin{align*}
H_{4k-1} &= \frac{H_{2k-1}}{2} + 2\sum_{1\le j\le k-1}
H_{2j}H_{4k-2j-2}
          +\alpha\cdot H_{4k-2} - \frac{{H_{2k-1}}^2}{2}
      -\sum_{1\le j\le 2k-2} H_j H_{4k-j-2} \\
H_{4k}   &= -\beta\cdot H_{4k-2} - \sum_{1\le j < 2k-1} H_j
H_{2(2k-j-1)}\\ H_{4k+1} &= \frac{H_{2k}}{2} + 2\sum_{1\le j <k}
H_{2j}H_{4k-2j}
          +\alpha\cdot H_{4k} + \frac{{H_{2k}}^2}{2}
      -\sum_{1\le j <2k} H_j H_{4k-j} \\
H_{4k+2} &= -\beta\cdot H_{4k} - \sum_{1\le j < 2k} H_j
H_{2(2k-j)}
\end{align*}
where $\alpha=-{\cal N}(j_{1,N})(\frac{1+N/2}{N})$ and
      $\beta =-{\cal N}(j_{1,N})(\frac{1}{N/2})$.
 Furthermore, we verified in \cite{Kim-KooB} that the above recursion
 can be also applied to 14 monstrous functions of even levels
 (including $\cal
 N(j_{1,2})$ and ${\cal N}(j_{1,6})$)
  which are Thompson series of type
 $2B$, $6C$, $6E$, $6F$, $10B$, $10E$, $14B$, $18C$, $18D$, $22B$,
 $30C$, $30G$, $42C$, $46AB$ (these are all replicable functions)
 and one monster-like function of type $18e$
 (for the definition of monster-like
 function, we refer to \cite{Ferenbaugh}).
 Therefore the Hauptmoduln mentioned above which have
 self-recursion formulas can be determined just by the first four
 coefficients $H_1$, $H_2$, $H_3$ and $H_4$ without the aid of
 2-plicate. What is more interesting would be the fact that there seems
 to be a connection between super-replicable functions and
 infinite dimensional Lie superalgebras. That is, considering the
 arguments from Borcherds \cite{Borcherds}, Kang \cite{Kang} and
 Koike \cite{Koike} we believe that the super-replication formulae
 in (\ref{super}) might suggest the existence of certain
 infinite dimensional Lie superalgebra whose denominator identity
 implies such formulae.

 \par Lastly, as an application of super-replication formulae
  we consider the following periodically
 vanishing property.
 Many of monstrous functions, for example, Thompson series of type
 $4B$, $4C$, $4D$, $6F$, $8B$, $8C$, $8D$, $8E$, $8F$, $9B$, etc
 have periodically vanishing properties among the Fourier coefficients
 (see the
 Table 1 in \cite{Mckay-Strauss}). This result must be known to
 experts, but we could not find a reference. Hereby we describe it
 in Theorem \ref{VVV}.
 Meanwhile, as for the case of super-replicable functions, we see from the
 Appendix, Table 4 that only the Haupmodul
 ${\cal N}(j_{1,12})$ seems to have such property.
 To this end, we shall first derive in \S 4
 an identity (\ref{twisted-replication}) which
 is analogous to the ``$2^k$-plication formula" (\cite{Ferenbaugh2}, \cite{Koike})
 satisfied
 by replicable functions. And, combining this with the super-replication formulae
  we are able to verify that the Fourier
 coefficients $H_m$ of ${\cal N}(j_{1,12})$ vanish whenever $m\equiv
 4 \mod 6$ (Corollary \ref{Replication3}).

 \par Through the article we adopt the following notations:
 \par\noindent $\bullet$ \, $S_{\Gamma_1(N)}$ \q the set of
 $\Gamma_1(N)$-inequivalent cusps
 \par\noindent $\bullet$ \, $q_h=e^{2\pi iz/h}, \q z\in {\frak H}$
 \par\noindent $\bullet$ \, $f|_{\sm a&b \\ c&d
 \esm}=f\left(\frac{az+b}{cz+d}\right)$
 \par\noindent $\bullet$ \,
     $f(z)=g(z)+O(1) \q \text{ means that $f(z)-g(z)$ is bounded as $z$
                    goes to $i\infty$. }$
\section{Hauptmodul of $\Gamma_1(12)$}
 \par
 In this section we investigate the generalities of the modular
 function $j_{1,12}$ which is under primary consideration and
 construct the Hauptmodul ${\cal N}(j_{1,12})$. We also examine
 some number theoretic property of ${\cal N}(j_{1,12})$. As for
 more arithmetic properties, we refer to \cite{Hong-Koo}.
\begin{lem}
Let $\frac ac$ and $\frac{a'}{c'}$ be fractions in lowest terms.
Then $\frac ac$ is $\Gamma_1(N)$-equivalent to $\frac{a'}{c'}$ if
and only if \q $\pm\sm a' \\ c' \esm \equiv
           \sm a+nc \\ c \esm \mod N$
           for some $n\in{\Bbb Z}$.
       \label{Equiv2}
\end{lem}
\begin{pf} Straightforward.
\end{pf}
Using the above lemma \, we can check that the cusps
$$ 0, \, 1/2,  \, 1/3, \, 1/4, \, 1/5, \, 1/6, \, 1/8, \,
   1/9, \, \infty, \, 5/12 $$
are $\Gamma_1(12)$-inequivalent. But from \cite{Chang2} we know
that the cardinality of $S_{\Gamma_1(12)}$ is 10, whence
$$ S_{\Gamma_1(12)}
   = \{0, \, 1/2, \, 1/3, \, 1/4, \, 1/5, \, 1/6, \, 1/8, \,
       1/9, \, \infty, \, 5/12 \}. $$
For later use we are in need of calculating the widths of the
cusps of $\Gamma_1(12)$.
\begin{lem} Let \, $a/c\in {\Bbb P}^1({\Bbb Q})$ be a cusp where $(a,c)=1$. Then
            the width of \, $a/c$ \, in $X_1(N)$ is given by $N/(c,N)$
        if $N\neq 4$.
        \label{Width}
\end{lem}
\begin{pf} If $N\mid 4$, the statement is obvious.
Hence, we assume that $N$ does not divide 4, i.e., $N\neq 1,2,4$.
First, choose $b$ and $d$ such that $\left(\begin{smallmatrix} a&
b \\ c& d \end{smallmatrix}\right)\in SL_2({\Bbb Z})$. Let $h$ be
the width of the cusp $a/c$. Then $h$ is the smallest positive
integer such that
$$ \m a&b \\ c&d \em
   \m 1&h \\ 0&1 \em
   {\m a&b \\ c&d \em}^{-1} \in \pm \Gamma_1(N). $$
Thus we have
$$ \m 1-cah & * \\ -c^2h & 1+cah \em \in \pm \Gamma_1(N). $$
If $\m 1-cah & * \\ -c^2h & 1+cah \em$ is an element of
$-\Gamma_1(N)$, by taking trace $2\equiv -2 \mod N$; hence $N\mid
4$. Thus when $N\neq 1,2,4$, $\m 1-cah & * \\ -c^2h & 1+cah \em
\in \Gamma_1(N)$. This condition is equivalent to saying that
$$ h\in \frac{N}{(c^2,N)} {\Bbb Z} \, \bigcap \,
        \frac{N}{(ca,N)} {\Bbb Z} \, = \, \frac{N}{(c,N)} {\Bbb Z}. $$
\end{pf}
We then have the following table of inequivalent cusps of
$\Gamma_1(12)$:
\begin{center}
{\bf Table 1}. \par
\begin{tabular}{|r|r|r|r|r|r|r|r|r|r|r|} \hline
 cusp & $\infty$ & $0$ & $1/2$ & $1/3$ & $1/4$ & $1/5$
         & $1/6$ & $1/8$ & $1/9$ & $5/12$ \\ \hline
 width &   1      &  12  &    6   &   4  &   3   &    12
          &   2  &    3   &   4  &   1    \\ \hline
\end{tabular}
\end{center}
Recall the Jacobi theta functions $\theta_2,\theta_3,$ and
$\theta_4$ defined by
\begin{align*}
    \theta_2(z) & =  \sum_{n\in{\Bbb Z}} q_2^{(n+\frac{1}{2})^2} \\
    \theta_3(z) & =  \sum_{n\in{\Bbb Z}} q_2^{n^2} \\
    \theta_4(z) & =  \sum_{n\in{\Bbb Z}} (-1)^nq_2^{n^2}
\end{align*}
   for $z\in{\frak H}$.
We have the following transformation formulas (\cite{Rankin}
pp.218-219).
\begin{align}
       \theta_2(z+1) &= e^{\frac{1}{4}\pi i}\theta_2(z) \label{T1} \\
       \theta_3(z+1) &= \theta_4(z) \label{T2} \\
       \theta_4(z+1) &= \theta_3(z)  \label{T3} \\
       \theta_2\left(-\frac{1}{z}\right)
       &= (-iz)^{\frac{1}{2}}\theta_4(z) \label{I1} \\
       \theta_3\left(-\frac{1}{z}\right)
       &= (-iz)^{\frac{1}{2}}\theta_3(z)  \label{I2} \\
       \theta_4\left(-\frac{1}{z}\right)
       &= (-iz)^{\frac{1}{2}}\theta_2(z). \label{I3}
\end{align}
\begin{lem} Let $k$ be an odd positive integer and
       $N$ be a multiple of $4$. Then for
       $F(z)\in M_{\frac k2}(\widetilde{\Gamma}_0(N))$
       and $m\ge 1$,
       $F(mz)\in M_{\frac k2}(\widetilde{\Gamma}_0(mN),\chi_m)$
       with $\chi_m(d)=\left(\frac{m}{d}\right)$ and $(d,m)=1$.
       \label{half2}
\end{lem}
\par\smallskip\noindent{\bf{\it Proof.}}\hskip0.4cm
\cite{Shimura2}, Proposition 1.3.
\par\noindent
Put $j_{1,12}(z)=\theta_3(2z)/\theta_3(6z)$.
\begin{theo}
(a) $\theta_3(2z)\in M_{\frac 12}(\widetilde{\Gamma}_0(4))$ and
    \, $\theta_3(6z)\in M_{\frac 12}(\widetilde{\Gamma}_0(12),\chi_3)$.
\par\noindent
(b) $K(X_1(12))$ is equal to \, ${\Bbb C}(j_{1,12}(z))$.
    $j_{1,12}$ takes
    the following value at each cusp:
    $j_{1,12}(\infty)=1$, $j_{1,12}(0)=\sqrt{3}$,
    $j_{1,12}(\frac12)=0$ (a simple zero),
    $j_{1,12}(\frac13)=i$, $j_{1,12}(\frac14)=\sqrt{3} i$,
    $j_{1,12}(\frac15)=-\sqrt3$,
    $j_{1,12}(\frac16)=\infty$ (a simple pole),
    $j_{1,12}(\frac18)=-\sqrt3 i$,
    $j_{1,12}(\frac19)=-i$,
    $j_{1,12}(\frac{5}{12})=-1$.
\label{Haupt2}
\end{theo}
\begin{pf} For the first part, we recall that (\cite{Koblitz}, p.184)
$$ \theta_3(2z)\in M_{\frac 12}(\widetilde{\Gamma}_0(4)).$$
Then by Lemma \ref{half2} we immediately get that
$$ \theta_3(6z)\in M_{\frac 12}(\widetilde{\Gamma}_0(12),\chi_3).$$
By the assertion (a), it is clear that $j_{1,12}(z)\in
K(X_1(12)).$ Thus for (b), it is enough to show that \,
$j_{1,12}(z)$ has only one simple zero and one simple pole on the
curve $X_1(12)$. As is well-known, $\theta_3(z)$ never vanishes
on ${\frak H}$. Hence we are forced to investigate the zeroes and
poles of \, $j_{1,12}$ at each cusp of $\Gamma_1(12)$. Let $S=\sm
0 & -1 \\ 1 & 0 \esm$ and $T=\sm 1 & 1 \\ 0 & 1 \esm$.
\par\noindent
(i) $s=\infty$:
   $$ j_{1,12}(\infty)=\lim_{z\to i\infty}\frac{\theta_3(2z)}{\theta_3(6z)}
              =\lim_{q\to 0}\frac{1+2q+2q^4+\cdots}{1+2q^3+2q^{12}+\cdots}
              =1. $$
\par\noindent
(ii) $s=0$:
 \begin{align*}
   j_{1,12}(0)  &=\lim_{z\to i\infty}
             {\frac{\theta_3(2z)}{\theta_3(6z)}\vline}_{ \q S} \\
           &=\lim_{z\to i\infty}
         \frac{\sqrt{-i\frac z2} \q \theta_3(\frac z2)}
                      {\sqrt{-i\frac z6} \q \theta_3(\frac z6)}
            \q \q \text{ by (\ref{I2})}  \\
               &=\sqrt{3}.
 \end{align*}
\par\noindent
(iii) $s=\frac12$:
 We Observe that $(ST^{-2}S)\infty=\sm -1&0\\-2&-1 \esm\infty=\frac12$.
 \par\noindent \q Considering the identities
 \begin{align*}
           \theta_3(2z)^2|_{[S]_1}
           &= z^{-1}\theta_3\left(-\frac{2}{z}\right)^2
           = z^{-1}\lbrace\left(-\frac{iz}{2}\right)
           ^{\frac{1}{2}}\theta_3\left(\frac{z}{2}\right)\rbrace^2
           \text{ \q by (\ref{I2})} \\
           &= -\frac i2 \theta_3\left(\frac{z}{2}\right)^2 \\
           \theta_3(2z)^2|_{[ST^{-2}]_1}
           &= -\frac{i}{2}\theta_3\left(\frac{z}{2}\right)^2|_{[T^{-2}]_1}
               = -\frac{i}{2}\theta_4\left(\frac{z}{2}\right)^2
               \text{ \q by (\ref{T3}) }\\
           \theta_3(2z)^2|_{[ST^{-2}S]_1}
           &= -\frac{i}{2}\theta_4\left(\frac{z}{2}\right)^2|_{[S]_1}
           = -\frac{i}{2}z^{-1}\{(-2iz)^{\frac{1}{2}}\theta_2(2z)\}^2
           \text{ \q by (\ref{I3})}\\
           &= -\theta_2(2z)^2,
     \end{align*}
     we get that
     \begin{align*}
           \theta_3(2z)^2|_{s=\frac12}
           &=\lim_{z\to i\infty}\theta_3(2z)^2|_{[ST^{-2}S]_1}
              =\lim_{z\to i\infty}-\theta_2(2z)^2   \\
           &=\lim_{z\to i\infty} -2^2 \q q_2(1+q^2+q^6+q^{12}+\cdots)^2
              \q \text{ since } \theta_2(z)= 2q_8(1+q+q^3+\cdots) \\
           &=0 \text{\q (a triple zero). }
     \end{align*}
     On the other hand
     \begin{align*}
           \theta_3(2z)^2|_{\sm -1&0\\-2&-1 \esm_1}
       & =\frac{1}{\sqrt3}\, \theta_3(2z)^2|_{\sm 3&0\\0&1 \esm_1
       \sm -1&0\\-2&-1 \esm_1}
         =\frac{1}{\sqrt3}\, \theta_3(2z)^2
              |_{\sm 1&1\\0&1 \esm_1\sm -1&0\\-2&-1 \esm_1
          \sm 1&-1\\0&3 \esm_1} \\
       & =\frac{1}{\sqrt3}\, \theta_3(2z)^2|_{\sm -1&0\\-2&-1 \esm_1
       \sm 1&-1\\0&3 \esm_1}
         =-\frac{1}{\sqrt3}\, \theta_2(2z)^2
              |_{\sm 1&-1\\0&3 \esm_1} \\
       & =-3^{-1} \, \theta_2\left(2\cdot \frac{z-1}{3}\right)^2
         =-3^{-1} e^{-\frac{\pi i}{3}} \, q_6(1+\cdots)^2,
     \end{align*}
 so that $\theta_3(6z)^2$ has a simple zero at $\frac12$. Thus $j_{1,12}^2
 =\frac{\theta_3(2z)^2}{\theta_3(6z)^2}$ has a double zero at $\frac12$, whence
 $j_{1,12}$ has a simple zero at $\frac12$.

\par\noindent
(iv) $s=\frac 13$: $(ST^{-3}S) \infty=\frac13$.
 First we recall that
  ${\frac{\theta_3(2z)}{\theta_3(6z)}\vline}_{ \q S}
   =\sqrt3 \, \frac{\theta_3(\frac z2)}
                      {\theta_3(\frac z6)}
            \q \q \text{ by (\ref{I2})}.  $
 Observe that $\theta_2(2z)=\frac12 (\theta_3(\frac z2)-\theta_4(\frac z2))$
             and $\theta_3(2z)=\frac12 (\theta_3(\frac z2)+\theta_4(\frac z2))$.
 From these identities we can write
      \, $\theta_3(\frac z2)=\theta_2(2z)+\theta_3(2z)$ and
      \, $\theta_3(\frac z6)=\theta_2(\frac23 z)+\theta_3(\frac23 z)$.
 Then we have that
 \begin{align*}
   {\frac{\theta_3(2z)}{\theta_3(6z)}\vline}_{ \q ST^{-3}}
    &= \sqrt{3}\cdot{\frac{\theta_2(2z)+\theta_3(2z)}{\theta_2(\frac23 z)+\theta_3(\frac23 z)}
               \vline}_{ \q T^{-3}}   \\
    &=\sqrt{3}\cdot\frac{(e^{-\frac{\pi i}{4}})^6 \theta_2(2z)+\theta_3(2z)}
                 {(e^{-\frac{\pi i}{4}})^2 \theta_2(\frac23 z)+\theta_3(\frac23 z)}
       \text{ \q by (\ref{T1}), (\ref{T2}) and (\ref{T3}),}\\
    &=\sqrt3\cdot\frac{i\theta_2(2z)+\theta_3(2z)}{-i\theta_2(\frac23 z)+\theta_3(\frac23 z)}
 \end{align*}
 and
 \begin{align*}
   {\frac{\theta_3(2z)}{\theta_3(6z)}\vline}_{ \q ST^{-3}S}
  &= \sqrt3\cdot
    {\frac{i\theta_2(2z)+\theta_3(2z)}{-i\theta_2(\frac23 z)+\theta_3(\frac23 z)}\vline}
    _{\q S} \\
  &= \sqrt3\cdot\frac{i\theta_4(\frac z2)+\theta_3(\frac z2)}
     {-i\sqrt3\theta_4(\frac32 z)+\sqrt3\theta_3(\frac32 z)}
     \text{ \q by (\ref{I1}) and (\ref{I2}) }
 \end{align*}
 which goes to $\frac{i+1}{-i+1}=i$ as $z\to i\infty$, so that
 $$ j_{1,12}\left(\frac13\right)=i.$$

\par\noindent
(v) $s=\frac14$: $\sm 1&0\\4&1 \esm \infty= \frac14$. In this
case we use the following well-known fact from \cite{Koblitz}
p.148 : For $\gamma\in\Gamma_0(4)$ and $z\in{\frak H}$,
$$ \Theta(\gamma z)=\left(\frac cd \right)\sqrt{\left(\frac{-1}{d}\right)}^{-1}
                    \sqrt{cz+d} $$
where $\Theta(z)=\theta_3(2z)$. Then
\begin{align*}
         {\frac{\theta_3(2z)}{\theta_3(6z)}\vline}_{ \q \sm 1&0\\4&1 \esm}
    &= \frac{\Theta(\sm 1&0\\4&1 \esm z)}{\Theta(\sm 3&0\\0&1 \esm \sm 1&0\\4&1 \esm z)}
     =\frac{\Theta(\sm 1&0\\4&1 \esm z)}{\Theta(\sm 3&2\\4&3 \esm \sm 1&-2\\0&3 \esm z)}\\
    &=\frac{\sqrt{4z+1} \q \Theta(z)}
       {\left(\frac43\right) i^{-1}\sqrt{4\cdot\frac{z-2}{3}+3}\q \Theta(\sm 1&-2\\0&3 \esm z)}\\
    &=\sqrt3 i \, \frac{\Theta(z)}{\Theta\left(\frac{z-2}{3}\right)}
\end{align*}
which tends to $\sqrt3 i$ when $z$ goes to $i\infty$. Therefore
 $$ j_{1,12}\left(\frac14\right)=\sqrt3 i.$$
\par\noindent
(vi) $s=\frac 15$: Because $\sm 5&1 \\ 24&5 \esm \sm 0&-1 \\ 1&0
\esm$ sends $\infty$ to $\frac 15$,
 \begin{align*}
   j_{1,12}\left(\frac15\right)  &=\lim_{z\to i\infty}
             {\frac{\theta_3(2z)}{\theta_3(6z)}\vline}_{
     \sm 5&1 \\ 24&5 \esm \sm 0&-1 \\ 1&0 \esm} \\
           &=\lim_{z\to i\infty}
         -{\frac{\theta_3(2z)}{\theta_3(6z)}\vline}_{\sm 0&-1\\1&0
               \esm}
            \q \q \text{ by (a)}  \\
               &=-j_{1,12}(0)=-\sqrt{3}.
 \end{align*}
(vii) $s=\frac16$: Observe that $(ST^{-6}S)\infty=\frac16$.
 \par\noindent \q Considering the identities
 \begin{align*}
   {\frac{\theta_3(2z)}{\theta_3(6z)}\vline}_{ \q S}
    &= \sqrt{3}\cdot\frac{\theta_3(\frac z2)}{\theta_3(\frac z6)}
       \text{ \q by } (\ref{I2})\\
   {\frac{\theta_3(2z)}{\theta_3(6z)}\vline}_{ \q ST^{-6}}
    &= \sqrt{3}\cdot{\frac{\theta_3(\frac z2)}{\theta_3(\frac z6)}
               \vline}_{ \q T^{-6}}
        =\sqrt{3}\cdot\frac{\theta_4(\frac z2)}{\theta_4(\frac z6)}
       \text{ \q by  (\ref{T2}) and (\ref{T3})} \\
   {\frac{\theta_3(2z)}{\theta_3(6z)}\vline}_{ \q ST^{-6}S}
    &= \sqrt{3}\cdot{\frac{\theta_4(\frac z2)}{\theta_4(\frac z6)}
               \vline}_{ \q S}
        =\sqrt{3}\cdot\frac{\sqrt{-2iz} \q \theta_2(2z)}
                   {\sqrt{-6iz} \q \theta_2(6z)}
       \text{ \q by } (\ref{I3}),
 \end{align*}
 we have that
 \begin{align*}
         j_{1,12}\left(\frac16\right)
     & =\lim_{z\to i\infty}
        {\frac{\theta_3(2z)}{\theta_3(6z)}\vline}_{ \q ST^{-6}S}
           =\lim_{z\to i\infty}
        \frac{\theta_2(2z)}{\theta_2(6z)} \\
     & =\lim_{z\to i\infty}
        \frac{2\,q_8^2(1+q^2+q^6+\cdots)}{2\,q_8^6(1+q^6+q^{18}+\cdots)}
        \q\q \text{ since $\theta_2(z)=2\, q_8(1+q+q^3+\cdots)$.}
 \end{align*}
 Thus by Table 1, $j_{1,12}$ has a simple pole at $\frac16$.
\par\noindent
(viii) $s=\frac 18$: Because $\sm -5&1 \\ -36&7 \esm \sm 1&0 \\
4&1 \esm$ sends $\infty$ to $\frac 18$,
 \begin{align*}
   j_{1,12}\left(\frac18\right)  &=\lim_{z\to i\infty}
             {\frac{\theta_3(2z)}{\theta_3(6z)}\vline}_{
     \sm -5&1 \\ -36&7 \esm \sm 1&0 \\ 4&1 \esm} \\
           &=\lim_{z\to i\infty}
         -{\frac{\theta_3(2z)}{\theta_3(6z)}\vline}_{\sm 1&0\\4&1
               \esm}
            \q \q \text{ by (a)}  \\
               &=-j_{1,12}\left(\frac14\right)=-\sqrt{3}i.
 \end{align*}
\par\noindent
(ix) $s=\frac 19$: Because $\sm -5&2 \\ -48&19 \esm \sm 1&0 \\
3&1 \esm$ sends $\infty$ to $\frac 19$,
 \begin{align*}
   j_{1,12}\left(\frac19\right)  &=\lim_{z\to i\infty}
             {\frac{\theta_3(2z)}{\theta_3(6z)}\vline}_{
     \sm -5&2 \\ -48&19 \esm \sm 1&0 \\ 3&1 \esm} \\
           &=\lim_{z\to i\infty}
         -{\frac{\theta_3(2z)}{\theta_3(6z)}\vline}_{\sm 1&0\\3&1
               \esm}
            \q \q \text{ by (a)}  \\
               &=-j_{1,12}\left(\frac13\right)=-i.
 \end{align*}
\par\noindent
(x) $s=\frac{5}{12}$: \q $\sm 5 &2 \\ 12&5 \esm \cdot
\infty=\frac{5}{12}$.
 \begin{align*}
   j_{1,12}\left(\frac{5}{12}\right)  &=\lim_{z\to i\infty}
             {\frac{\theta_3(2z)}{\theta_3(6z)}\vline}_{
     \sm 5&2 \\ 12&5  \esm} \\
           &=\lim_{z\to i\infty}
         -\frac{\theta_3(2z)}{\theta_3(6z)}
            \q \q \text{ by (a)}  \\
               &=-1.
 \end{align*}
\end{pf}
We will now construct the Hauptmodul ${\cal N}(j_{1,12})$ from
the modular function $j_{1,12}$ mentioned in Theorem \ref{Haupt2}.
\begin{align*}
\frac{2}{j_{1,12}(z)-1}
  & =\frac{2\q\theta_3(6z)}{\theta_3(2z)-\theta_3(6z)}
   =\frac{2(1+2q^3+2q^{12}+2q^{27}+\cdots)}
      {2q-2q^3+2q^4+2q^9-2q^{12}+\cdots} \\
  & =\frac1q+q+q^2+q^3-q^6-q^7-q^8-q^9+q^{11}+2q^{12}+\cdots,
\end{align*}
which is in \, $q^{-1}{\Bbb Z}[[q]]$. From the uniqueness of the
normalized generator it follows that ${\cal
N}(j_{1,12})=\frac{2}{j_{1,12}-1}$.
 By Theorem \ref{Haupt2}-(b) we have the following
table:
\begin{center}
{\bf Table 2}. Cusp values of $j_{1,12}$ and ${\cal N}(j_{1,12})$
\par
\begin{tabular}{|c|c|c|c|c|c|c|c|c|c|c|} \hline
 $s$ & $\infty$ & $0$ & $1/2$ & $1/3$ & $1/4$ & $1/5$
                      & $1/6$ & $1/8$ & $1/9$ & $5/12$  \\ \hline
 $j_{1,12}(s)$ & $1$ & $\sqrt3$ & $0$ & $i$ & $\sqrt3 i$ & $-\sqrt3$
                      & $\infty$ & $-\sqrt3 i$ & $-i$ & $-1$ \\ \hline
${\cal N}(j_{1,12})(s)$ & $\infty$ & $\sqrt3+1$ & $-2$ & $-1-i$ &
$\frac{-1-\sqrt3 i}{2}$ & $1-\sqrt3$
                      & $0$ & $\frac{-1+\sqrt3 i}{2}$ & $-1+i$ & $-1$ \\ \hline
\end{tabular}
\end{center}
\vskip0.5cm
\begin{theo}
       Let $d$ be a square free positive integer and \, $t=N(j_{1,N})$
       be the normalized generator of $K(X_1(N))$.
       Let $s$ be a cusp of \, $\Gamma_1(N)$
       whose width is $h_s$. If $t\in q^{-1}\Bbb Z[[q]]$ and
       $\prod_{s\in S_{\Gamma_1(N)}
        \smallsetminus \{\infty\}} (t(z)-t(s))^{h_s} $ is a polynomial in
       $\Bbb Z[t]$, then $t(\tau)$ is an algebraic integer
       for \, $\tau\in{\Bbb Q}(\sqrt{-d})\cap\frak H$.
       \label{Integer}
\end{theo}
\begin{pf} Let $j(z)=\dfrac1q+744+196884q+\cdots$, the elliptic
 modular function. It is well-known that $j(\tau)$ is an algebraic
 integer for $\tau\in{\Bbb Q}(\sqrt{-d})\cap\frak H$ (\cite{Lang},
 \cite{Shimura1}). For algebraic proofs, see \cite{Deuring},
 \cite{Neron}, \cite{Serre} and \cite{Silverman}. Now, we view $j$
 as a function on the modular curve $X_1(N)$. Then $j$ has a pole
 of order $h_s$ at the cusp $s$. On the other hand, $t(z)-t(s)$
 has a simple zero at $s$. Thus $$ j\times \prod_{s\in
 S_{\Gamma_1(N)}
      \smallsetminus \{\infty\}} (t(z)-t(s))^{h_s} $$
 has a pole only at $\infty$ whose degree is \,
 $\mu_N=[\overline{\Gamma}(1):
 \overline{\Gamma}_1(N)]$,
 and so by the Riemann-Roch Theorem it is a monic polynomial in \,
 $t$ \, of degree \, $\mu_N$ which we denote by $f(t)$. Since
 $\prod_{s\in S_{\Gamma_1(N)}
      \smallsetminus \{\infty\}} (t(z)-t(s))^{h_s} $
      is a polynomial in $\Bbb Z[t]$ and
$j$, $t$ have integer coefficients in the $q$-expansions, $f(t)$
is a monic polynomial in ${\Bbb Z}[t]$ of degree $\mu_N$. This
shows that $t(\tau)$ is integral over ${\Bbb Z}[j(\tau)]$.
Therefore $t(\tau)$ is integral over $\Bbb Z$
 \, for $\tau\in{\Bbb Q}(\sqrt{-d})\cap\frak H$.
\end{pf}
\begin{cor} For \, $\tau\in{\Bbb Q}(\sqrt{-d})\cap\frak H$,
${\cal N}(j_{1,12})(\tau)$ is an algebraic integer.
\end{cor}
\begin{pf}
${\cal N}(j_{1,12})$ has integral Fourier coefficients. And by
Table 1 and 2,
$$ \prod_{s\in S_{\Gamma_1(12)
      \smallsetminus \{\infty\}}} (t(z)-t(s))^{h_s}
   = (t^2-2t-2)^{12} \, (t+2)^6 \, (t^2+2t+2)^4 \, (t^2+t+1)^3 \,
     t^2 \, (t+1) \q \in \Bbb Z[t].  $$
Now the assertion is immediate from Theorem \ref{Integer}.
\end{pf}

\section{Super-replication formulae}
\par
Let $\Delta^n$ be the set of $2\times 2$ integral matrices $\sm
a&b
\\ c&d \esm$ where $a\in 1+N{\Bbb Z}, \, c\in N{\Bbb Z}$, and
$ad-bc=n$. Then $\Delta^n$ has the following right coset
decomposition: (See \cite{Koblitz}, \cite{Miyake},
\cite{Shimura1})
\begin{equation}
\Delta^n=\bigcup_{\smallmatrix a|n \\ (a,N)=1 \endsmallmatrix}
     \bigcup_{i=0}^{\frac{n}{a}-1} \Gamma_1(N)\sigma_a
          \m a&i \\ 0&\frac{n}{a} \em   \label{Coset}
\end{equation}
where $\sigma_a\in SL_2({\Bbb Z})$ such that $\sigma_a\equiv \sm
a^{-1}&0 \\ 0&a \esm \mod N$. Let $f(z)$ be a modular function
with respect to $\Gamma_1(N)$. For brevity, let us call it {\it
$f(z)$ is on} $\Gamma_1(N)$.  For $f\in K(X_1(N))$ we define an
operator $U_n$ and $T_n$ by
\begin{align*}
 f|_{U_n} &=n^{-1} \sum_{i=0}^{n-1} f|_{\sm 1&i\\0&n \esm} \\
 \intertext{and}
 f|_{T_n} &= n^{-1} \sum_{\smallmatrix a\mid n \\ (a,N)=1
\endsmallmatrix}
   \sum_{i=0}^{\frac na-1} f|_{\sigma_a \sm a& i \\ 0& \frac na \esm}.
\end{align*}
\begin{lem}
For $f\in K(X_1(N))$ and $\gamma_0\in \Gamma_0(N)$,
$(f|_{T_n})|_{\gamma_0}=(f|_{\gamma_0})|_{T_n}$ for any positive integer $n$.
In particular, $f|_{T_n}$ is again on $\Gamma_1(N)$.
\label{Hecke}
\end{lem}
\begin{pf}
First we claim that
$$ \Delta^n\gamma_0=\gamma_0\Delta^n \text{ \q
 for }\gamma_0=\sm x&y\\z&w\esm\in\Gamma_0(N).$$
Let $\sm a&b\\c&d \esm\in \Delta^n$. Then
 $$
{\gamma_0}^{-1}\sm a&b\\c&d \esm \gamma_0 = \sm w& -y\\ -z&x\esm
\sm a&b\\c&d \esm \sm x&y\\z&w\esm \equiv \sm 1&*\\ 0&d \esm \mod
N.
 $$
Hence \, ${\gamma_0}^{-1}\Delta^n\gamma_0 \subset \Delta^n$ \, so
that \, $\Delta^n\gamma_0\subset \gamma_0\Delta^n$. By the same
argument we can show the reverse inclusion. We note that
\begin{align*}
\Delta^n\gamma_0 & =
     \bigcup_{\smallmatrix a|n \\ (a,N)=1 \endsmallmatrix}
     \bigcup_{i=0}^{\frac{n}{a}-1} \Gamma_1(N)\sigma_a
          \m a&i \\ 0&\frac{n}{a} \em  \gamma_0 \\
\intertext{ and }
\gamma_0\Delta^n & =
     \bigcup_{\smallmatrix a|n \\ (a,N)=1 \endsmallmatrix}
     \bigcup_{i=0}^{\frac{n}{a}-1} \gamma_0\Gamma_1(N)\sigma_a
          \m a&i \\ 0&\frac{n}{a} \em  \\
 & = \bigcup_{\smallmatrix a|n \\ (a,N)=1 \endsmallmatrix}
     \bigcup_{i=0}^{\frac{n}{a}-1} \Gamma_1(N)\gamma_0\sigma_a
          \m a&i \\ 0&\frac{n}{a} \em
      \text{ \q because $\Gamma_1(N)\triangleleft\Gamma_0(N)$.}
\end{align*}
Here we note that $\sigma_a\sm a&i \\ 0&\frac{n}{a} \esm
\gamma_0$'s are the matrices appearing in the definition of
$(f|_{T_n})|_{\gamma_0}$ and $\gamma_0\sigma_a\sm a&i \\
0&\frac{n}{a} \esm$'s are those appearing in the definition of
$(f|_{\gamma_0})|_{T_n}$. Now the assertion follows.
\end{pf}
For a positive integer $N$ with $g_{1,N}=0$, we let $t$ (resp.
$t_0$) be the Hauptmodul of $\Gamma_1(N)$ (resp. $\Gamma_0(N)$).
And, we write $X_n(t)=\frac 1n q^{-n} + \sum_{m\ge 1} H_{m,n} q^m$
and $X_n(t_0)=\frac 1n q^{-n} + \sum_{m\ge 1} h_{m,n} q^m$.
\begin{lem}
 For positive integers $m$ and $n$, $H_{m,n}=H_{n,m}$ and
 $h_{m,n}=h_{n,m}$. \label{CCC}
\end{lem}
\begin{pf}
 Let $p=e^{2\pi i y}$ and $q=e^{2\pi i z}$ with $y,z\in {\frak
 H}$. Note that $X_n(t)$ can be viewed as the coefficient of
 $p^n$-term in $-\log p - \log (t(y)-t(z))$ (\cite{Norton}). Thus
 $H_{m,n}$ becomes the coefficient of $p^n q^m$-term of
 \begin{align*}
  & -\log p - \log (t(y)-t(z)) \\
  = & -\log (1-p/q)+\log (p^{-1} - q^{-1}) - \log (t(y)-t(z)) \\
  = & \sum_{i\ge 1} \frac 1i (p/q)^i - F(p,q)
 \end{align*}
 where $F(p,q)=\log \left(
  \frac{p^{-1}-q^{-1}+\sum_{i\ge 1} H_i (p^i -q^i)}{p^{-1}-q^{-1}}
  \right)$.
 We then come up with $F(p,q)=F(q,p)$, which implies that
 $H_{m,n}=H_{n,m}$. Similarly if we
 work with $t_0$ instead of $t$, the identity $h_{m,n}=h_{n,m}$ follows.
 \end{pf}
\begin{theo}
For positive integers $n$ and $l$ such that $(n,N)=(l,n)=1$,
$$ X_l(t)|_{T_n} = X_{ln}(t)|_{\sigma_n} + c $$
where $c$ is a constant. In particular,
$$ t|_{T_n} = X_{n}(t)|_{\sigma_n} + c. $$
\label{Hecke2}
\end{theo}
\begin{pf}
Since $X_l(t)$ has poles only at $\Gamma_1(N)\infty$, the poles
of $X_l(t)|_{T_n}$ can occur only at $\sm a&i \\ 0&\frac na
\esm^{-1} {\sigma_a}^{-1} \Gamma_1(N)\infty$ where $a$ and $i$ are
the indices appearing in the definition of $T_n$. On the other
hand, we have
\begin{align*}
\sm a&i \\ 0&\frac na \esm^{-1} {\sigma_a}^{-1}\Gamma_1(N)\infty
  &=n^{-1}\sm \frac na&-i \\ 0&a \esm{\sigma_a}^{-1}\Gamma_1(N)\infty \\
  & = \sm \frac na&-i \\ 0&a \esm{\sigma_a}^{-1}\Gamma_1(N)\infty.
\end{align*}
Let $\gamma$ be an element in $\Gamma_1(N)$.
Then
\begin{align*}
\sm \frac na&-i \\ 0&a \esm{\sigma_a}^{-1}\gamma\infty
& \equiv \sm \frac na&-i \\ 0&a \esm
  \sm a&0\\0&a^{-1}\esm\sm 1&*\\0&1\esm\infty \mod N \\
& \equiv \sm n&*\\0&1 \esm \infty = \frac{n+Nm}{Nk}
\end{align*}
for some $k,m\in {\Bbb Z}$. If there exists an integer $x>1$ with
$x\mid (n+Nm,Nk)$, then $x$ must divide $\det\left[\sm \frac
na&-i\\ 0&a \esm {\sigma_a}^{-1}\gamma\right]=n$. In this case
$x\nmid N$ because $(n,N)=1$. Therefore $\sm a&i \\ 0&\frac na
\esm^{-1} {\sigma_a}^{-1}\gamma\infty$ is of the form
$\gamma_0\infty$ for some $\gamma_0\in \Gamma_0(N)$. Now we
conclude that $X_l(t)|_{T_n}$ can have poles only at
$\gamma_0\infty$ for some $\gamma_0\in\Gamma_0(N)$. By Lemma
\ref{Hecke},
\begin{align*}
 & n(X_l(t)|_{T_n})|_{\gamma_0} = n(X_l(t)|_{\gamma_0})|_{T_n} \\
= & \sum_{a\mid n }
   \sum_{i=0}^{\frac na-1} (X_l(t)|_{\gamma_0})
   |_{\sigma_a \sm a& i \\ 0& \frac na \esm}
= \sum_{a\mid n }
   \frac na {(X_l(t)|_{\gamma_0})
   |_{\sigma_a}}|_{U_\frac na} (az) \\
=& \sum_{\smallmatrix a\mid n \\ a\neq n \endsmallmatrix}
   \frac na {(X_l(t)|_{\gamma_0})
   |_{\sigma_a}}|_{U_\frac na} (az) + X_l(t)|_{\gamma_0\sigma_n}(nz).
\end{align*}
Here we note that
$$ \sum_{\smallmatrix a\mid n \\ a\neq n \endsmallmatrix}
   \frac na {(X_l(t)|_{\gamma_0})
   |_{\sigma_a}}|_{U_\frac na} (az) = O(1).$$
 In fact, if $\gamma_0\sigma_a \notin\pm\Gamma_1(N)$, it is clear
 that ${(X_l(t)|_{\gamma_0})
   |_{\sigma_a}}|_{U_\frac na}$ has a holomorphic $q$-expansion.
   Otherwise
$$ X_l(t)|_{\gamma_0\sigma_a U_\frac na} = X_l(t)|_{U_\frac na} =
(l^{-1}q^{-l}+\text{ terms of positive degree })|_{U_\frac na} =
O(1)
$$ because $(l,n)=1$. Now we have
\begin{equation}
 n(X_l(t)|_{T_n})|_{\gamma_0}
= X_l(t)|_{\gamma_0\sigma_n}(nz) + O(1).
\label{A}
\end{equation}
This implies that $X_l(t)|_{T_n}$ has a pole at $\gamma_0\infty$
if and only if $\gamma_0\sigma_n \in \pm\Gamma_1(N)$, that is,
$\gamma_0\in \pm\Gamma_1(N){\sigma_n}^{-1}$. Hence $X_l(t)|_{T_n}$
has poles only at cusps $\Gamma_1(N){\sigma_n}^{-1}\infty$. In
this case we derive from (\ref{A})
\begin{equation*}
  n(X_l(t)|_{T_n})|_{{\sigma_n}^{-1}}
= X_l(t)|_{{\sigma_n}^{-1}\sigma_n}(nz)+ O(1) = l^{-1}q^{-ln} + O(1).
\end{equation*}
Thus \, $(X_l(t)|_{T_n})|_{{\sigma_n}^{-1}} =(nl)^{-1}q^{-ln} +
O(1)$ \, and \, $(X_l(t)|_{T_n})|_{{\sigma_n}^{-1}}$ \, has poles
only at cusps $\sigma_n
\Gamma_1(N){\sigma_n}^{-1}\infty=\Gamma_1(N)\infty.$ On the other
hand, $X_{ln}(t)$ has poles only at $\Gamma_1(N)\infty$ too and
$X_{ln}(t)=(ln)^{-1}q^{-ln} + O(1).$ Therefore
$(X_l(t)|_{T_n})|_{{\sigma_n}^{-1}} = X_{ln}(t)+c$ for some
constant $c$. Then we have $X_l(t)|_{T_n} =
X_{ln}(t)|_{\sigma_n}+c$, as desired.
\end{pf}
\begin{cor} Let $N$ be a positive interger such that the genus
$g_{1,N}$ is zero and
$[\overline{\Gamma}_0(N):\overline{\Gamma}_1(N)]\le 2$. For
positive integers $n,l$ and $m$ such that $(n,N)=(l,n)=1$, we have
$$ \sum_{\smallmatrix e\mid (m,n)\\ e>0 \endsmallmatrix}
   e^{-1}\left\{ \psi (e)
      \left( 2 H_{\frac{mn}{e^2},l}-h_{\frac{mn}{e^2},l}\right)
       +h_{\frac{mn}{e^2},l}\right\}
      = \psi (n)(2H_{m,ln}-h_{m,ln})+h_{m,ln} $$
where $\psi: \, ({\Bbb Z}/N{\Bbb Z})^\times \to \{ \pm 1\}$ is a
character defined by
$$ \psi(e)= \begin{cases} 1, & \text{ if } e\equiv \pm 1 \mod N \\
-1, & \text{otherwise.} \end{cases} $$ \label{RRR1}
\end{cor}
\begin{pf}
It follows from Theorem \ref{Hecke2} that
\begin{equation}
 X_l(t)|_{T_n} = \sum_{\smallmatrix e\mid n\\ e>0 \endsmallmatrix}
   e^{-1} \left(X_l(t)|_{\sigma_e}\right)|_{U_\frac ne}(ez)
   = X_{ln}(t)|_{\sigma_n}+ \text{ constant. } \label{BB}
\end{equation}
 Note that for each positive integer r,
$$ X_r(t)|_{\sigma_e}+X_r(t)=
\begin{cases} 2 X_r(t), &  \text{ if } e\equiv \pm 1 \mod N \\
               X_r(t_0)+ \text{ constant }, &  \text{ otherwise.}
\end{cases} $$
In the above when $e$ is not congruent to $\pm 1 \mod N$, $
X_r(t)|_{\sigma_e}+X_r(t)$ is on $\Gamma_0(N)$ and has poles only
at $\Gamma_0(N)\infty$ with $r^{-1}q^{-r}$ as its pole part, which
guarantees the above equality. We then have $$
X_r(t)|_{\sigma_e}=\frac12
  \{\psi(e)(2X_r(t)-X_r(t_0))+X_r(t_0)\} + \text{ constant.} $$
Now (\ref{BB}) reads
\begin{gather*}
\sum_{\smallmatrix e\mid n \\ e>0 \endsmallmatrix}
 e^{-1}\cdot \frac 12
  \{\psi(e)(2X_l(t)-X_l(t_0))+X_l(t_0)\}|_{U_\frac ne}(ez) \\
 = \frac 12
  \{\psi(n)(2X_{ln}(t)-X_{ln}(t_0))+X_{ln}(t_0)\}+ \text{constant.}
\end{gather*}
Comparing the coefficients of $q^m$-terms on both sides, we get
the corollary.
\end{pf}
\begin{cor}
 Let $N$ be a positive integer such that the genus $g_{1,N}$ is zero
 and $[\overline{\Gamma}_0(N)
 :\overline{\Gamma}_1(N)]\le 2$. For positive integers $a,b,c,d$ with
 $ab=cd$, $(a,b)=(c,d)$ and $(b,N)=(d,N)=1$,
 \begin{equation}
 H_{a,b}=\psi(bd) H_{c,d} + \frac{(1-\psi(bd))}{2}h_{c,d}.
 \label{super}
 \end{equation}
 \label{Repli2}
\end{cor}
\begin{pf}
In Corollary \ref{RRR1} we take $n=b$, $l=1$ and $m=a$. Then it
follows from the conditions and the replicability of $h_{m,n}$
that
$$ \psi(b)(2H_{a,b}-h_{a,b})=\psi(d)(2H_{c,d}-h_{c,d}). $$
Now the assertion follows.
\end{pf}
\begin{cor}
Let $N$ be a positive integer with $g_{1,N}=0$ and
$[\overline{\Gamma}_0(N):\overline{\Gamma}_1(N)]=2$. If $(mn,N)=1$
and $mn\not\equiv \pm 1 \mod N$, then $h_{m,n}=2H_{m,n}$.
\end{cor}
\begin{pf}
In Corollary \ref{Repli2} we take $a=m, b=n$ and $c=n, d=m$.
 The condition that $\psi(mn)=-1$ implies
\begin{align*}
 H_{m,n}  &= -H_{n,m}+h_{n,m} \\
          &= -H_{m,n}+h_{m,n} \text{ \q by Lemma \ref{CCC}. }
\end{align*}
This proves the corollary.
\end{pf}

\section{Vanishing property in the Fourier coefficients of ${\cal N}(j_{1,12})$}
\par
 As mentioned in the introduction, many of the Thompson series
 have periodically
 vanishing properties among the Fourier coefficients. Now we will
 give a more theoretical explanation for these phenomena.
 \par Let $T_g$ be the Thompson series
 of type $g$ and $\Gamma_g$ be its corresponding genus zero group.
 To describe $\Gamma_g$ we are in need of some notations.
 Let $N$ be a positive integer and $Q$ be any Hall divisor of
$N$, that is, $Q$ be a positive divisor of $N$ for which $(Q,
N/Q)=1$. We denote by $W_{Q,N}$ a matrix $\sm Qx & y \\ Nz & Qw
\esm$ with $\det W_{Q,N}=Q$ and $x,y,z$ and $w\in \Bbb Z$, and
call it an Atkin-Lehner involution. Let $S$ be a subset of Hall
divisors of $N$ and let $\Gamma=N+S$ be the subgroup of
$PSL_2(\Bbb R)$ generated by $\Gamma_0(N)$ ($=\{ \sm a&b\\c&d \esm
\in SL_2(\Bbb Z) \, | \, c\equiv 0 \mod N \}$) and all
Atkin-Lehner involutions $W_{Q,N}$ for $Q\in S$.
 For a positive divisor $h$ of $24$, let
$n$ be a multiple of $h$ and set $N=nh$. When $S$ is a subset of
Hall divisors of $n/h$, we denote by $\Gamma_0(n|h)+S$ the group
generated by $\sm h&0 \\ 0&1 \esm^{-1} \circ \Gamma_0(n/h)
\circ\sm h&0 \\ 0&1 \esm $ and $\sm h&0 \\ 0&1 \esm^{-1} \circ
W_{Q,n/h} \circ \sm h&0 \\ 0&1 \esm$ for all $Q\in S$. If there
exists a homomorphism $\lambda$ of $\Gamma_0(n|h)+S$ into $\Bbb
C^*$ such that
 \begin{align}
 & \lambda(\Gamma_0(N))=1, \label{H1} \\
 &\lambda(\sm 1&1/h\\ 0&1 \esm)=e^{-2\pi i/h}, \label{H2} \\
 & \lambda(\sm 1&0\\ n&1 \esm)= \begin{cases} e^{2\pi i/h} & \text{ \q if $n/h\in
   S$}, \\ e^{-2\pi i/h} & \text{ \q if $n/h \notin S$,}
     \end{cases} \label{H3} \\
 & \lambda \text{ is trivial on all Atkin-Lehner involutions of
    } \Gamma_0(N) \text{ in } \Gamma_0(n|h)+S, \label{H4}
 \end{align}
  then we let $n|h+S$ be the kernel of $\lambda$ which is
 a subgroup of $\Gamma_0(n|h)+S$ of index $h$.
 Ferenbaugh (\cite{Ferenbaugh}) found out a necessary and
 sufficient condition for the homomorphism $\lambda$ to exist and
 calculated the genera of groups of type $n|h+S$. All the genus
 zero groups of type $n|h+S$ are listed in \cite{Ferenbaugh},
 Table 1.1 and 1.2.
 Now we have the
 following theorem.
 \begin{theo} Write $T_g(q)=q^{-1}+\sum_{m\ge 1} c_g(m) q^m$.
 \par\noindent
 (i) If $\Gamma_g=\Gamma_0(N)$ and $h$ is the largest
 integer such that $h\mid 24$ and $h^2 \mid N$, then
 $ c_g(m)=0 \text{ unless } m\equiv -1 \mod h. $
 \par\noindent
 (ii) If $\Gamma_g=N+S$ and there exists a prime $p$ such that
 $p^2|N$ and $p\nmid Q$ for all $Q\in S$, then
 $ c_g(m)=0 \text{ whenever } m\equiv 0 \mod p. $
 \par\noindent
 (iii) If $\Gamma_g=n|h+S$, then
 $ c_g(m)=0 \text{ unless } m\equiv -1 \mod h. $
 \label{VVV}
 \end{theo}
 \begin{pf} (i) Note that $\sm 1& h^{-1}\\ 0&1 \esm $ belongs to
 the normalizer of $\Gamma_0(N)$. Thus $T_g|_{\sm 1& h^{-1}\\ 0&1
 \esm}$ has poles only at $\infty$, which is a simple pole.
 This enables us to write $T_g|_{\sm 1& h^{-1}\\ 0&1
 \esm}=c\cdot T_g$ for some constant $c$. By comparing the
 coefficients of $q^m$-terms on both sides, the assertion follows.
 \par\noindent
 (ii) From Corollary 3.1 \cite{Koike} it follows that $T_g|_{U_p}=0$.
 Thus
 (ii) is clear.
 \par\noindent
 (iii) Considering the identity in \cite{Koike}, p.27 we have
 $T_g(z+1/h)=e^{-2\pi i/h}\cdot T_g(z)$. Then
 $T_g(q)=q^{-1}+\sum_{0<l\in {\Bbb Z}} c_g(lh-1) q^{lh-1}$, which
 implies (iii).
 \end{pf}
 Unlike the cases of Thompson series, when we handle the
 super-replicable function ${\cal N}(j_{1,12})$ we can not
 directly use the ingredients adopted in Theorem \ref{VVV}.
 Therefore we start with
 \begin{lem} For $\sm a&b\\c&d \esm\in \Gamma_0(12)$,
 $j_{1,12}|_{\sm a&b\\c&d \esm}
 =\left(\frac{3}{d}\right)j_{1,12}$. Here $\left(\frac{\cdot}{\cdot}\right)$
 denotes the generalized quadratic residue symbol.
 \label{Haupt}
 \end{lem}
 \begin{pf} Immediate from Lemma \ref{half2} and Theorem
 \ref{Haupt2}.
 \end{pf}
 We fix $N=12$ and let $t$ denote the Hauptmodul ${\cal N}(j_{1,12})$
 in what follows.

\begin{lem} For $\sm a&b \\ c&d \esm \in\Gamma_0(24)$,
$$ \left(t|_{U_2}\right)_\chi|_{\sm a&b\\c&d \esm}
= (-1)^\frac{c}{24}\cdot \left(\frac{3}{a+\frac c4}\right)
  \cdot\left(t|_{U_2}\right)_\chi$$
where $\chi=\left(\frac{-1}{\cdot}\right)$ is the Jacobi symbol
and $\left(t|_{U_2}\right)_\chi$ is the twist of $t|_{U_{2}}$ by
$\chi$ (\cite{Koblitz}, p.127). \label{AAA}
\end{lem}
\begin{pf}
From \cite{Koblitz}, p.128 we observe that
\begin{align}
\left(t|_{U_2}\right)_\chi &=\frac{1}{\sqrt{-4}}
\left(t|_{U_2}\left(z+\frac14\right)-t|_{U_2}\left(z+\frac34\right)\right)
\label{DD} \\
&=\frac{1}{\sqrt{-4}} \sum_{i\in ({\Bbb Z}/4{\Bbb Z})^\times}
\left(\frac{-1}{i}\right) \left(t|_{U_2}+\frac12\right)|_{\sm
4&i\\0&4 \esm}. \notag
\end{align}
It then follows from \cite{Kim-KooA}, Corollary 28 \, that
$$ \left(t-t\left(\frac{1}{6}\right)\right)\times t|_{U_2} =H_2.
$$ If we compare the coefficients of $q$-term on both sides, we get
$H_4-t\left(\frac{1}{6}\right)\cdot H_2 = 0$. And, substituting
$H_2=1$ and $H_4=0$ we get
 \, $t\left(\frac{1}{6}\right) = 0$.
Now
\begin{equation}
t|_{U_2}=\frac1t=\frac{j_{1,12}-1}{2}.
\label{KK}
\end{equation}
 Then for $\sm a&b\\c&d \esm\in\Gamma_0(12)$
\begin{align}
\left(t|_{U_2}+\frac12\right)|_{\sm a&b\\c&d \esm}
&= \frac12 j_{1,12}|_{\sm a&b\\c&d \esm}
 =\frac12\cdot\left(\frac{3}{d}\right)j_{1,12}
  \text{ \q by Lemma \ref{Haupt}}
  \label{EE} \\
&=\left(\frac{3}{d}\right)\left(t|_{U_2}+\frac12\right)
 =\left(\frac{3}{a}\right)\left(t|_{U_2}+\frac12\right) \notag
\end{align} since $ad\equiv 1 \mod 12$.
Let $\sm a&b\\c&d \esm \in\Gamma_0(24)$. For $i=1,3$, we consider
$\sm 4&i\\0&4 \esm\sm a&b\\c&d \esm = \sm 4a+ic&4b+id\\4c&4d
\esm$. Since $(4a+ic,4c)$ divides $\det\sm 4a+ic&4b+id\\4c&4d
\esm$, we must have $(4a+ic,4c)=4$; hence $(a+ic/4,c)=1$. Thus we
can choose integers $x_i$ and $y_i$ such that $\gamma_i=\sm
a+ic/4&x_i\\c&y_i \esm\in SL_2({\Bbb Z})$. Write
${\gamma_i}^{-1}\sm4&i\\0&4 \esm\sm a&b\\c&d \esm = \sm
4&z_i\\0&4\esm$ for some integer $z_i$. Then we have
$$ \sm4&i\\0&4 \esm\sm a&b\\c&d \esm
 = \sm a+ic/4&x_i\\c&y_i \esm\sm 4&z_i\\0&4\esm$$
for some integer $z_i$. Comparing (1,2)-component on both sides we get
$$ 4b+id=(a+ic/4)z_i + 4x_i.$$ Thus $id\equiv(a+ic/4)z_i \mod 4$. Then
\begin{align*}
z_i & \equiv (a+ic/4)id \mod 4 \text{ \q \q because \q $n^2\equiv
1 \mod 4$
       \q for every odd integer $n$} \\
    & \equiv aid+i^2\cdot\frac c4\cdot d \mod 4 \equiv i + \frac c4\cdot d
       \mod 4 \text{ \q \q due to \, $ad\equiv 1 \mod 4$} \\
    & \equiv i+6c_1d \mod 4, \text{\q \q where we write $c=24c_1$.}
\end{align*}
Therefore we derive that for $\sm a&b\\c&d \esm \in\Gamma_0(24)$,
\begin{align*}
\sqrt{-4}\cdot\left(t|_{U_2}\right)_\chi|_{\sm a&b\\c&d\esm} &=
\sum_{i\in ({\Bbb Z}/4{\Bbb Z})^\times} \left(\frac{-1}{i}\right)
\left(t|_{U_2}+\frac12\right)|_{\sm 4&i\\0&4 \esm\sm a&b\\c&d
\esm}
\text{ \q by (\ref{DD})} \\
&= \sum_{i\in ({\Bbb Z}/4{\Bbb Z})^\times}
\left(\frac{-1}{i}\right) \left(t|_{U_2}+\frac12\right)|_
{\sm a+ic/4&x_i\\c&y_i \esm\sm 4&z_i\\0&4 \esm} \\
&= \sum_{i\in ({\Bbb Z}/4{\Bbb Z})^\times}
\left(\frac{-1}{i}\right)\left(\frac{3}{a+ic/4}\right)
\left(t|_{U_2}+\frac12\right)|_{\sm 4&z_i\\0&4 \esm} \text{ \q by
(\ref{EE}).}
\end{align*}
Now if we set $c=24c_1$ as before, then we have
\begin{align*}
& \q \sqrt{-4}\cdot\left(t|_{U_2}\right)_\chi|_{\sm a&b\\c&d\esm}
\\ &= \sum_{i\in ({\Bbb Z}/4{\Bbb Z})^\times}
\left(\frac{-1}{i}\right)\left(\frac{3}{a+6c_1i}\right)
\left(t|_{U_2}+\frac12\right)|_{\sm 4&i+6c_1d\\0&4 \esm}  \\ &=
\left(\frac{3}{a+6c_1}\right)\sum_{i\in ({\Bbb Z}/4{\Bbb
Z})^\times} \left(\frac{-1}{i}\right)
\left(t|_{U_2}+\frac12\right)|_{\sm 4&i+6c_1d\\0&4 \esm}  \\ &
\text{ \q \q since $a+6c_1i\equiv a+6c_1 \mod 12$} \\ &=
\left(\frac{3}{a+6c_1}\right)\sum_{i\in ({\Bbb Z}/4{\Bbb
Z})^\times}
\left(\frac{-1}{i}\right)\left(\frac{-1}{i+6c_1d}\right)
\left(\frac{-1}{i+6c_1d}\right)
\left(t|_{U_2}+\frac12\right)|_{\sm 4&i+6c_1d\\0&4 \esm}  \\ &=
(-1)^{c_1}\cdot\left(\frac{3}{a+6c_1}\right)\cdot\sqrt{-4}\cdot
\left(t|_{U_2}\right)_\chi
\end{align*}
because $\left(\frac{-1}{i}\right)\left(\frac{-1}{i+6c_1d}\right)
=\left(\frac{-1}{i^2+6c_1id}\right)=\left(\frac{-1}{1+2c_1id}\right)
=(-1)^{c_1id}=(-1)^{c_1}$ and $i+6c_1d$ runs over $({\Bbb Z}/4\Bbb
Z)^\times$. This completes the lemma.
\end{pf}
\begin{lem}
 (i) For each $k\ge 1$, we put $g(z)=(-1)^{k-1}
 \cdot \frac 12\cdot \left(X_{2^k}(t)(z)-X_{2^k}(t)\left(z+\frac 12\right)\right).$
  Then $g$ belongs to $K(X_1(24))$.
 \par\noindent
 (ii) For $\sm a&b \\ c&d \esm \in\Gamma_0(24)$ and $k\ge 1$, $$
 \left(t|_{U_{2^k}}\right)_\chi|_{\sm a&b\\c&d \esm} =
 (-1)^\frac{c}{24}\cdot \left(\frac{3}{a+\frac c4}\right)
  \cdot\left(t|_{U_{2^k}}\right)_\chi.$$ In particular,
  $\left(t|_{U_{2^k}}\right)_\chi$ lies in $K(X_1(24))$.
 \label{Replication5}
\end{lem}
\begin{pf}
First we note that for $n\mid N^\infty$, $T_n=U_n$. Here, by
$n\mid N^\infty$ we mean that $n$ divides some power of $N$. To
show $g\in K(X_1(24))$, we observe that
$$ g=(-1)^{k-1}
\left(X_{2^k}(t)(z)-X_{2^k}(t)|_{U_2}(2z)\right).$$ Then using
Lemma \ref{Hecke} we obtain that $g\in K(X_1(24))$. For $\sm a&b
\\ c&d \esm \in\Gamma_0(12)$,
\begin{align}
\left(t|_{U_{2^k}}+\frac12\right)|_{\sm a&b\\c&d \esm} &
=\left(t|_{U_{2}}+\frac12\right)|_{U_{2^{k-1}}\sm a&b\\c&d \esm}
 =\left(t|_{U_{2}}+\frac12\right)|_{\sm a&b\\c&d \esm
U_{2^{k-1}}}
   \text{ \q by Lemma \ref{Hecke}}
   \label{FF} \\
& =\left(\frac 3a\right)\left(t|_{U_{2}}+\frac12\right)|_{U_{2^{k-1}}}
   \text{ \q by (\ref{EE})} \notag \\
& =\left(\frac 3a\right)\left(t|_{U_{2^k}}+\frac12\right). \notag
\end{align}
Now we can proceed in the same manner as in the proof of Lemma
\ref{AAA}.
\end{pf}
\begin{lem}
For $k\ge 1$,
\begin{align*}
& (i) \q \left(t|_{U_{2^k}}\right)|_{\sm 1&0\\6&1 \esm } =
\begin{cases}
\frac12 {q_2}^{-1}+O(1), & \text{ if $k=1$} \\
O(1), & \text{ otherwise. }
\end{cases}
\\
& (ii) \q \left(t|_{U_{2^k}}\right)|_{\sm 1&0\\3&1 \esm } =
\begin{cases}
O(1), & \text{ if $k=1$} \\
\frac14 {q_4}^{-1}+O(1), & \text{ if $k=2$ } \\
-\frac i8 {q_2}^{-1}+O(1), & \text{ if $k=3$ } \\
(-1)^{k-1}\cdot\frac{i}{2^k}\cdot{q}^{-2^{k-4}}+O(1), & \text{ if $k\ge 4$. }
\end{cases}
\end{align*}
\label{Replication6}
\end{lem}
\begin{pf}
(i) First, $t|_{\sm 1&0\\6&1 \esm}$ is holomorphic at $\infty$
because $t$ has poles only at the cusps $\Gamma_1(12)\infty$. Now
for $k\ge 1$,
\begin{align*}
\left(t|_{U_{2^k}}\right)|_{\sm 1&0\\6&1 \esm }
&=\left(\left(t|_{U_{2^{k-1}}}\right)|_{U_2}\right)|_{\sm 1&0\\6&1 \esm }
=\frac12\left(t|_{U_{2^{k-1}}}\right)|_{\sm 1&0\\0&2 \esm \sm 1&0\\6&1 \esm}
+\frac12\left(t|_{U_{2^{k-1}}}\right)|_{\sm 1&1\\0&2 \esm \sm 1&0\\6&1 \esm}
\\
&=\frac12\left(t|_{U_{2^{k-1}}}\right)|_{\sm 1&0\\12&1 \esm \sm 1&0\\0&2 \esm}
+\frac12\left(t|_{U_{2^{k-1}}}\right)|_{\sm 7&4\\12&7 \esm \sm 1&-1\\0&2 \esm}
\\
&=\frac12\left(t|_{U_{2^{k-1}}}\right)\left(\frac z2\right)
+\frac12\left(\left(t|_{\sm 7&4\\12&7 \esm}\right)|_
{U_{2^{k-1}}}\right)|_{\sm 1&-1\\0&2 \esm}
\text{ \q by Lemma \ref{Hecke} } \\
&=\begin{cases}
\frac12 {q_2}^{-1}+O(1), & \text{ if $k=1$} \\
O(1), & \text{ otherwise. }
\end{cases}
\end{align*}
\par\noindent
(ii) We observe that $t|_{\sm 1&0\\3&1 \esm}\in O(1)$. And for
$k\ge 1$,
\begin{align}
\left(t|_{U_{2^k}}\right)|_{\sm 1&0\\3&1 \esm }
&=\left(\left(t|_{U_{2^{k-1}}}\right)|_{U_2}\right)|_{\sm 1&0\\3&1 \esm }
\label{HH} \\
& =\frac12\left(t|_{U_{2^{k-1}}}\right)|_{\sm 1&0\\0&2 \esm \sm 1&0\\3&1 \esm}
+\frac12\left(t|_{U_{2^{k-1}}}\right)|_{\sm 1&1\\0&2 \esm \sm 1&0\\3&1 \esm}
\notag \\
&=\frac12\left(t|_{U_{2^{k-1}}}\right)|_{\sm 1&0\\6&1 \esm \sm 1&0\\0&2 \esm}
+\frac12\left(t|_{U_{2^{k-1}}}\right)|_{\sm 2&1\\3&2 \esm \sm 2&0\\0&1 \esm}
\notag
\end{align}
has a holomorphic Fourier expansion if $k=1$. Thus we suppose
$k\ge 2$. We then derive that
\begin{align*}
\left(t|_{U_{2^{k-1}}}\right)|_{\sm 2&1\\3&2 \esm }
&=\left(\left(t|_{U_{2^{k-2}}}\right)|_{U_2}\right)|_{\sm 2&1\\3&2 \esm } \\
&=\frac12\left(t|_{U_{2^{k-2}}}\right)|_{\sm 1&0\\0&2 \esm \sm 2&1\\3&2 \esm}
+\frac12\left(t|_{U_{2^{k-2}}}\right)|_{\sm 1&1\\0&2 \esm \sm 2&1\\3&2 \esm}
\\
&=\frac12\left(t|_{U_{2^{k-2}}}\right)|_{\sm 1&0\\3&1 \esm \sm 2&1\\0&1 \esm}
+\frac12\left(t|_{U_{2^{k-2}}}\right)|_
{\sm 11&-1\\12&-1 \esm \sm 1&0\\6&1 \esm\sm 1&1\\0&2 \esm}
\\
&=\frac12\left(t|_{U_{2^{k-2}}}\right)|_{\sm 1&0\\3&1 \esm \sm 2&1\\0&1 \esm}
+\frac12\left(t|_{U_{2^{k-2}}}\right)|_{\sm 1&0\\6&1 \esm \sm 1&1\\0&2 \esm}.
\end{align*}
If we substitute the above into (\ref{HH}), for $k\ge 2$,
\begin{align}
 & \q  \left(t|_{U_{2^{k}}}\right)|_{\sm 1&0\\3&1 \esm } \label{II}\\
 & =\frac12\left(t|_{U_{2^{k-1}}}\right)|_{\sm 1&0\\6&1 \esm
   \sm 1&0\\0&2 \esm }
  +\frac14\left(t|_{U_{2^{k-2}}}\right)|_{\sm 1&0\\6&1 \esm \sm 1&1\\0&2 \esm \sm 2&0\\0&1 \esm }
  +\frac14\left(t|_{U_{2^{k-2}}}\right)|_{\sm 1&0\\3&1 \esm \sm 2&1\\0&1 \esm \sm 2&0\\0&1 \esm }
   \notag \\
 & =\frac12\left(t|_{U_{2^{k-1}}}\right)|_{\sm 1&0\\6&1 \esm
   \sm 1&0\\0&2 \esm }
 +\frac14\left(t|_{U_{2^{k-2}}}\right)|_{\sm 1&0\\6&1 \esm \sm 2&1\\0&2 \esm }
 +\frac14\left(t|_{U_{2^{k-2}}}\right)|_{\sm 1&0\\3&1 \esm \sm 4&1\\0&1 \esm }
\notag
\end{align}
When $k=2$,
\begin{align*}
& \q  \left(t|_{U_{4}}\right)|_{\sm 1&0\\3&1 \esm }
=\frac12\left(t|_{U_{2}}\right)|_{\sm 1&0\\6&1 \esm}\left(\frac
z2\right)
 +\frac14 t|_{\sm 1&0\\6&1 \esm }\left(z+\frac12\right)
 +\frac14 t|_{\sm 1&0\\3&1 \esm }\left(4z+1\right)  \\
& = \frac14 {q_4}^{-1}+O(1) \text{ \q\q by (i). }
\end{align*}
If $k=3$,
\begin{align*}
& \q  \left(t|_{U_{8}}\right)|_{\sm 1&0\\3&1 \esm }
=\frac12\left(t|_{U_{4}}\right)|_{\sm 1&0\\6&1 \esm}\left(\frac
z2\right)
 +\frac14\left(t|_{U_{2}}\right)|_{\sm 1&0\\6&1 \esm } \left(z+\frac12\right)
 +\frac14\left(t|_{U_{2}}\right)|_{\sm 1&0\\3&1 \esm }(4z+1) \\
& = \frac18 e^{-\pi i\left(z+\frac12\right)} + O(1) =
-\frac{i}{8}{q_2}^{-1}+O(1) \text{ \q\q by (i) and the case $k=1$
in (ii). }
\end{align*}
For $k\ge 4$, we will show by induction on $k$ that
\begin{equation}
\left(t|_{U_{2^{k}}}\right)|_{\sm 1&0\\3&1 \esm }
   = (-1)^{k-1}\cdot \frac{i}{2^k}\cdot q^{-2^{k-4}}+O(1).
   \label{JJ}
\end{equation}
First we note that by (i) and (\ref{II})
$$\left(t|_{U_{2^{k}}}\right)|_{\sm 1&0\\3&1 \esm }
  =\frac14\left(t|_{U_{2^{k-2}}}\right)|_{\sm 1&0\\3&1 \esm }(4z+1)+O(1). $$
If $k=4$,
\begin{align*}
\left(t|_{U_{2^{4}}}\right)|_{\sm 1&0\\3&1 \esm }
  &=\frac14\left(t|_{U_{2^{4-2}}}\right)|_{\sm 1&0\\3&1 \esm }(4z+1)+O(1)\\
&=\frac{1}{16}e^{-\frac{\pi i}{2}(4z+1)}+O(1)
  \text{ \q \q by the case $k=2$} \\
&= (-1)^{4-1}\cdot \frac{i}{2^4}\cdot q^{-2^{4-4}}+O(1).
\end{align*}
Thus when $k=4$, (\ref{JJ}) holds. Meanwhile, if $k=5$ then we get
that
\begin{align*}
\left(t|_{U_{2^{5}}}\right)|_{\sm 1&0\\3&1 \esm }
  &=\frac14\left(t|_{U_{2^{5-2}}}\right)|_{\sm 1&0\\3&1 \esm }(4z+1)+O(1)\\
&=-\frac{i}{32}e^{-\pi i(4z+1)}+O(1)
  \text{ \q \q by the case $k=3$} \\
&= (-1)^{5-1}\cdot \frac{i}{2^5}\cdot q^{-2^{5-4}}+O(1).
\end{align*}
Therefore in this case (\ref{JJ}) is also valid. Now for $k\ge 6$,
\begin{align*}
\left(t|_{U_{2^{k}}}\right)|_{\sm 1&0\\3&1 \esm }
  &=\frac14\left(t|_{U_{2^{k-2}}}\right)|_{\sm 1&0\\3&1 \esm }(4z+1)+O(1)\\
&=\frac{1}{4}\cdot (-1)^{k-2-1}\cdot \frac{i}{2^{k-2}}\cdot
 e^{-2\pi i(4z+1)\cdot 2^{k-2-4}}+O(1) \\
& \text{ \q \q \, by induction hypothesis for $k-2$} \\
&= (-1)^{k-1}\cdot \frac{i}{2^k}\cdot q^{-2^{k-4}}+O(1).
\end{align*}
This proves the lemma.
\end{pf}
\begin{theo} For $k\ge 1$,
\begin{equation}
 \left(t|_{U_{2^k}}\right)_\chi(z)
=(-1)^{k-1}\cdot \frac12\cdot \left(
X_{2^k}(t)(z)-X_{2^k}(t)\left(z+\frac12\right)\right).
\label{twisted-replication}
\end{equation}
This identity twisted by a character $\chi$ is analogous to the
``$2^k$-plication formula" (\cite{Ferenbaugh2}, \cite{Koike})
satisfied by the Hauptmodul of $\Gamma_0(N)$. \label{CC}
\end{theo}
\begin{pf}
Put $g(z)= (-1)^{k-1}\cdot \frac12\cdot \left(
X_{2^k}(t)(z)-X_{2^k}(t)\left(z+\frac12\right)\right)$ as before.
 We see by Lemma \ref{Replication5} that both $\left(t|_{U_{2^k}}\right)_\chi$
 and $g$ sit in $K(X_1(24))$.
 For $\left(t|_{U_{2^k}}\right)_\chi = g$, we will show that
 $\left(t|_{U_{2^k}}\right)_\chi - g$ has no poles in ${\frak H}^*$.
 Recall that $$ \left(t|_{U_{2^k}}\right)_\chi
 =\frac{1}{\sqrt{-4}} \sum_{i\in ({\Bbb Z}/4{\Bbb Z})^\times}
\chi(i)\left(t|_{U_{2^k}}\right)|_{\sm 4&i\\0&4 \esm}
 =\frac{1}{\sqrt{-4}}
\sum_{i\in ({\Bbb Z}/4{\Bbb Z})^\times}\sum_{j=0}^{2^k-1}
\chi(i)t|_{\sm 1&j\\0&2^k \esm\sm 4&i\\0&4 \esm}.
$$
Since $t$ has poles only at $\Gamma_1(12)\infty$,
$\left(t|_{U_{2^k}}\right)_\chi$ can also have poles only at
\begin{align*}
\sm 4&i\\0&4 \esm^{-1}\sm 1&j\\0&2^k \esm^{-1}\Gamma_1(12)\infty
& = 16^{-1}\cdot 2^{-k}\cdot
  \sm 4&-i\\0&4 \esm\sm 2^k&-j\\0&1 \esm\Gamma_1(12)\infty \\
& =
  \sm 2^{k+2}&-4j-i\\0&4 \esm\Gamma_1(12)\infty.
\end{align*}
Let $\sm a&b\\c&d \esm\in \Gamma_1(12)$. Then
$$ \sm 2^{k+2}&-4j-i\\0&4 \esm\sm a&b\\c&d \esm\infty
   =\sm 2^{k+2}a-4jc-ic & * \\ 4c&* \esm \infty
   =\frac{2^{k+2}a-4jc-ic}{4c}. $$
Observe that $(2^{k+2}a-4jc-ic,4c)\mid 2^k \cdot 16=2^{k+4}$.
Write $(2^{k+2}a-4jc-ic,4c)=2^l$ for some integer $l\ge 0$. Then
   $$s=\frac{(2^{k+2}a-4jc-ic)/2^l}{4c/2^l} $$
is in lowest terms. Since $12\mid c$, $s$ is of the form
$s=\frac{n}{3m}$ for some integers $m$ and $n$. We assume
$(3m,n)=1$. Here we consider two cases.
\par\noindent
(i) $2^2\nmid m$: \par Choose integers $x$ and $y$ such that $\sm
n&x\\ 3m&y\esm\in SL_2({\Bbb Z})$ and consider
\begin{align*}
\left(t|_{U_{2^k}}\right)_\chi|_{\sm n&x\\3m&y \esm }
&=\frac{1}{\sqrt{-4}}
\left(t|_{U_{2^k}\sm 4&1\\0&4 \esm\sm n&x\\3m&y \esm}
     -t|_{U_{2^k}\sm 4&3\\0&4 \esm\sm n&x\\3m&y \esm}\right) \\
&=\frac{1}{\sqrt{-4}}
\left(t|_{U_{2^k}\sm 4n+3m&*\\12m&* \esm}
     -t|_{U_{2^k}\sm 4n+9m&*\\12m&* \esm}\right).
\end{align*}
Since $(3m,n)=1$ and $2^2\nmid m$, we can write $ \sm
4n+3m&*\\12m&* \esm=\gamma_0 U_1$ and $ \sm 4n+9m&*\\12m&*
\esm={\gamma_0}' U_2$ where both $\gamma_0$ and ${\gamma_0}'$ are
in $\Gamma_0(12)$ and $U_1$, $U_2$ are upper triangular matrices.
Then by (\ref{FF}), $\left(t|_{U_{2^k}}\right)_\chi|_{\sm
n&x\\3m&y \esm} \in O(1)$. Hence, if $2^2\nmid m$ then
$\left(t|_{U_{2^k}}\right)_\chi$ is holomorphic at the cusp
$s=\frac{n}{3m}$.
\par\noindent
(ii) $2^2\mid m$: \par If $2^3\mid m$, then $s$ is of the form
$s=\sm n&x\\3m&y \esm\infty$ with $\sm
n&x\\3m&y\esm\in\Gamma_0(24)$. Thus by Lemma \ref{Replication5},
$$ \left(t|_{U_{2^k}}\right)_\chi|_{\sm n&x\\3m&y \esm}
= (-1)^\frac{3m}{24}\cdot \left(\frac{3}{n+3m/4}\right)
  \cdot\left(t|_{U_{2^k}}\right)_\chi\in O(1).$$
As for the other cases, if we use Lemma \ref{Equiv2}, it is easy
to see that $s$ is equivalent to $\frac{1}{12}$ or $\frac{5}{12}$
under $\Gamma_1(24)$.
\par
Thus we conclude that $\left(t|_{U_{2^k}}\right)_\chi$ can have
poles only at $\frac{1}{12}$, $\frac{5}{12}$ under
$\Gamma_1(24)$-equivalence. Next, let us investigate the poles of
$g$. Recall that $X_{2^k}(t)$ has poles only at
$\Gamma_1(12)\infty$. Therefore $g$ can have poles only at $\sm
2&i\\0&2 \esm^{-1}\Gamma_1(12)\infty$ for $i=0,1$. And, for $\sm
a&b\\c&d \esm\in\Gamma_1(12)$,
\begin{align*}
\sm 2&i\\0&2 \esm^{-1}\sm a&b\\c&d \esm\infty & = \sm 2&-i\\0&2
\esm\sm a&b\\c&d \esm\infty
 = \frac{2a-ic}{2c}
 = \frac{a-ic/2}{c} \text{ \q \q in lowest terms. }
\end{align*}
Hence by Lemma \ref{Equiv2}, $g$ can have poles only at
$\frac{1}{24}$, $\frac{5}{24}$, $\frac{7}{24}$, $\frac{11}{24}$,
$\frac{1}{12}$, $\frac{5}{12}$ under $\Gamma_1(24)$-equivalence.
At $\frac{1}{24}$, $\frac{5}{24}$, $\frac{7}{24}$ and
$\frac{11}{24}$, it is easy to check that $g$ is holomorphic. For
example, at $\frac{5}{24}$,
\begin{align*}
g|_{\sm 5&1\\24&5 \esm}
& = (-1)^{k-1}\cdot \frac12\cdot
\left( X_{2^k}(t)|_{\sm 5&1\\24&5 \esm}
  -X_{2^k}(t)|_{\sm 2&1\\0&2 \esm\sm 5&1\\24&5 \esm}\right) \\
& = (-1)^{k-1}\cdot \frac12\cdot
\left( X_{2^k}(t)|_{\sm 5&1\\24&5 \esm}
  -X_{2^k}(t)|_{\sm 34&7\\48&10 \esm}\right) \\
& = (-1)^{k-1}\cdot \frac12\cdot
\left( X_{2^k}(t)|_{\sm 5&1\\24&5 \esm}
  -X_{2^k}(t)|_{\sm 17&12\\24&17 \esm\sm 2&-1\\0&2 \esm}\right) \\
&  \in O(1) \text{ \q since $\frac{5}{24},
  \frac{17}{24}\notin\Gamma_1(12)\infty$. }
\end{align*}
Now it remains to show that $\left(t|_{U_{2^k}}\right)_\chi-g$ has
no poles at the cusps equivalent to $\frac{1}{12},\frac{5}{12}$
under $\Gamma_1(24)$. At $\frac{1}{12}$,
\begin{align}
\left(t|_{U_{2^k}}\right)_\chi|_{\sm 1&0\\12&1 \esm }
&=\frac{1}{\sqrt{-4}}
\left(t|_{U_{2^k}\sm 4&1\\0&4 \esm\sm 1&0\\12&1 \esm}
     -t|_{U_{2^k}\sm 4&3\\0&4 \esm\sm 1&0\\12&1 \esm}\right)
     \label{GG} \\
&=\frac{1}{\sqrt{-4}}
\left(t|_{U_{2^k}\sm 1&0\\3&1 \esm\sm 16&1\\0&1 \esm}
     -t|_{U_{2^k}\sm 5&-1\\6&-1 \esm\sm 8&1\\0&2 \esm}\right) \notag \\
&=\frac{1}{\sqrt{-4}}
\left(t|_{U_{2^k}\sm 1&0\\3&1 \esm\sm 16&1\\0&1 \esm}
     -t|_{U_{2^k}\sm 11&-1\\12&-1 \esm\sm 1&0\\6&1 \esm
                 \sm 8&1\\0&2 \esm}\right)
\notag \\
&=\frac{1}{\sqrt{-4}}
\left(t|_{U_{2^k}\sm 1&0\\3&1 \esm}(16z+1)
     -t|_{U_{2^k}\sm 1&0\\6&1 \esm}\left(4z+\frac12\right)\right).
\notag
\end{align}
By (\ref{GG}) and Lemma \ref{Replication6}, we have the following:
\begin{align*}
\text{\noindent If $k=1$, } \hskip2cm
\left(t|_{U_{2}}\right)_\chi|_{\sm 1&0\\12&1 \esm }
& = \frac{1}{\sqrt{-4}}\cdot
   \left(-\frac12 e^{-\pi i\left(4z+\frac12\right)}+O(1)\right) = \frac{1}{4}q^{-2}+O(1).
\end{align*}
\begin{align*}
\text{\noindent If $k=2$, } \hskip2cm
\left(t|_{U_{2^2}}\right)_\chi|_{\sm 1&0\\12&1 \esm }
& = \frac{1}{\sqrt{-4}}\cdot
   \left(\frac14 e^{-\frac{\pi i}{2}\left(16z+1\right)}+O(1)\right)
  = -\frac{1}{8}q^{-4}+O(1).
\end{align*}
\begin{align*}
\text{\noindent If $k=3$, } \hskip2cm
\left(t|_{U_{2^3}}\right)_\chi|_{\sm 1&0\\12&1 \esm }
& = \frac{1}{\sqrt{-4}}\cdot
   \left(-\frac i8 e^{-\pi i\left(16z+1\right)}+O(1)\right) = \frac{1}{16}q^{-8}+O(1).
\end{align*}
\begin{align*}
\text{\noindent If $k\ge 4$, } \hskip2cm
\left(t|_{U_{2^k}}\right)_\chi|_{\sm 1&0\\12&1 \esm }
& = \frac{1}{\sqrt{-4}}\cdot
   \left((-1)^{k-1}\cdot \frac{i}{2^k}\cdot
   e^{-2\pi i\cdot 2^{k-4}\cdot \left(16z+1\right)}+O(1)\right) \\
&=  (-1)^{k-1}\cdot \frac{1}{2^{k+1}}\cdot q^{-2^k}+O(1).
\end{align*}
Observe that the identities for $k=1,2$ and $3$ are the same as
the last one when $k\ge 4$. Hence we conclude that for all $k\ge
1$,
$$\left(t|_{U_{2^k}}\right)_\chi|_{\sm 1&0\\12&1 \esm }
 = (-1)^{k-1}\cdot \frac{1}{2^{k+1}}\cdot q^{-2^k}+O(1).$$
On the other hand,
\begin{align*}
g|_{\sm 1&0\\12&1 \esm}
& = (-1)^{k-1}\cdot \frac12\cdot
\left( X_{2^k}(t)|_{\sm 1&0\\12&1 \esm}
  -X_{2^k}(t)|_{\sm 2&1\\0&2 \esm\sm 1&0\\12&1 \esm}\right) \\
& = (-1)^{k-1}\cdot \frac12\cdot
\left( X_{2^k}(t)
  -X_{2^k}(t)|_{\sm 7&4\\12&7 \esm\sm 2&-1\\0&2 \esm}\right) \\
& = (-1)^{k-1}\cdot \frac{1}{2^{k+1}}\cdot q^{-2^k}+O(1).
\end{align*}
Thus
$$\left(\left(t|_{U_{2^k}}\right)_\chi-g\right)|_
  {\sm 1&0\\12&1 \esm }\in O(1).$$
At $\frac{5}{12}$, we see that $\sm 5&2\\12&5 \esm=\sm
-19&2\\-48&5 \esm\sm 1&0\\12&1 \esm$ sends $\infty$ to
$\frac{5}{12}$.  Then
\begin{align*}
 \left(t|_{U_{2^k}}\right)_\chi|_{\sm 5&2\\12&5 \esm}
&= \left(t|_{U_{2^k}}\right)_\chi|_
  {\sm -19&2\\-48&5 \esm\sm 1&0\\12&1 \esm} \\
&= (-1)^\frac{-48}{24}\cdot \left(\frac{3}{-19+\frac{-48}{4}}\right)
  \cdot\left(t|_{U_{2^k}}\right)_\chi|_{\sm 1&0\\12&1 \esm}
  \text{\q by Lemma \ref{Replication5}} \\
&= (-1)\cdot \left(t|_{U_{2^k}}\right)_\chi|_{\sm 1&0\\12&1 \esm }
 = (-1)^{k}\cdot \frac{1}{2^{k+1}}\cdot q^{-2^k}+O(1).
\end{align*} And
\begin{align*}
g|_{\sm 5&2\\12&5 \esm}
& = (-1)^{k-1}\cdot \frac12\cdot
\left( X_{2^k}(t)|_{\sm 5&2\\12&5 \esm}
  -X_{2^k}(t)|_{\sm 2&1\\0&2 \esm\sm 5&2\\12&5 \esm}\right) \\
& = (-1)^{k-1}\cdot \frac12\cdot
\left( X_{2^k}(t)|_{\sm 5&2\\12&5 \esm}
  -X_{2^k}(t)|_{\sm 11&-1\\12&-1 \esm\sm 2&1\\0&2 \esm}\right) \\
& = (-1)^{k}\cdot \frac{1}{2}\cdot X_{2^k}(t)|_{\sm 2&1\\0&2 \esm}+O(1)
= (-1)^{k}\cdot \frac{1}{2^{k+1}}\cdot q^{-2^k}+O(1).
\end{align*}
This implies that
$$\left(\left(t|_{U_{2^k}}\right)_\chi-g\right)|_
  {\sm 5&2\\12&5 \esm }\in O(1),$$
 from which the theorem follows.
\end{pf}
Now, we are ready to show periodically vanishing property of
${\cal N}(j_{1,12})$.
\begin{cor}
As before we let $t$ be the Hauptmodul of $\Gamma_1(12)$ and
write $X_n(t)=\frac 1n q^{-n}+\sum_{m\ge 1}H_{m,n} q^m$. Then
we have \par (i)
$ H_{m,2^k} = (-1)^{k-1}\left(\frac{-1}{m}\right)H_{2^km,1}
\text{ \q \q for odd $m$.} $
\par (ii)
$H_m=0$ whenever $m\equiv 4 \mod 6$, and $m=5$.
\label{Replication3}
\end{cor}
\begin{pf}
First if we compare the coefficients of $q^m$-terms on both sides
of the identity in Theorem \ref{CC}, we get (i).
 We see from the Appendix, Table 4 that $H_5=0$. On the other
hand, by the super-replication formula (Corollary \ref{Repli2})
it follows that for $m$ relatively prime to 12,
$$ H_{m,2^k}=H_{2^k,m}=
\begin{cases} H_{2^km,1}, & \text{ if } m\equiv \pm 1 \mod 12 \\
-H_{2^km,1}, & \text{ if } m\equiv \pm 5 \mod 12
\end{cases} $$
because $h_{2^km,1}=0$ in this case (\cite{Koike}, Corollary 3.1).
Then $H_{2^km,1}=0$ when $k$ is odd and $m\equiv 5 \mod 6$, or
$k$ is even and $m\equiv 1 \mod 6$. It is easy to see that
\begin{align*}
& \{ 2^km \, | \, k, m\ge 1, \,
                \text{$k$ odd, $m\equiv 5 \mod 6$} \} \cup
 \{ 2^km \, | \, k, m\ge 1, \,
                \text{$k$ even, $m\equiv 1 \mod 6$} \}  \\
=& \{ l\in{\Bbb Z} \, | \, l\equiv 4 \mod 6 \}.
\end{align*}
This proves (ii).
\end{pf}

 \vskip1cm \centerline{\bf Appendix. \q Fourier coefficients of the
 Hauptmodul ${\cal N} (j_{1,N})$}
\par We shall make use of the following modular forms to construct $j_{1,N}$.
 For $z\in {\frak H}$,
\par\noindent$\bullet$
$\eta(z)$ \q \q the Dedekind eta function
\par\noindent$\bullet$
$G_2(z)$ \q \q Eisenstein series of weight 2
\par\noindent$\bullet$
$G_2^{(p)}(z)=G_2(z)-pG_2(pz)$ for each prime $p$
\par\noindent$\bullet$
$E_2(z)= G_2(z)/(2\zeta(2))$ \q \q normalized Eisenstein series
of weight 2
\par\noindent$\bullet$
$E_2^{(p)}(z)=E_2(z)-pE_2(pz)$ for each prime $p$
\par\noindent$\bullet$
${\cal P}_{N,{\bf a}}(z)=\cal
P_{L_z}\left(\frac{a_1z+a_2}{N}\right) $ \q \q $N$-th division
value of ${\cal P}$ where ${\bf a}=(a_1,a_2), L_z={\Bbb Z} z+{\Bbb
Z}$ and ${\cal P}_L(\tau)$ is the Weierstrass ${\cal P}$-function
(relative to a lattice $L$)
\par Now we get the following tables due to \cite{Kim-Koo2}-\cite{Kim-Koo4}:
\vskip0.5cm
\begin{center}
 {\bf Table 3}. Hauptmoduln $\cal N (j_{1,N})$
                       \vskip0.2cm
\begin{tabular}
{|c||c|c|}\hline $N$ & $j_{1,N}$ & $\cal N (j_{1,N})$ \\
\hline\hline 1 & $j(z)$ & $j(z)-744$ \\ \hline 2   &
$\theta_2(z)^8/\theta_4(2z)^8$ & $256/j_{1,2}+24$ \\ \hline 3   &
$E_4(z)/E_4(3z)$ & $240/(j_{1,3}-1)+9$  \\ \hline 4   &
$\theta_2(2z)^4/\theta_3(2z)^4$ & $16/j_{1,4}-8$ \\ \hline 5   &
$\frac{4\eta(z)^5/\eta(5z)+E_2^{(5)}(z)}{\eta(5z)^5/\eta(z)}$ &
      $-8/(j_{1,5}+44)-5$  \\ \hline
6   &
$\frac{G_2^{(2)}(z)-G_2^{(2)}(3z)}{2G_2^{(2)}(z)-G_2^{(3)}(z)}$ &
      $2/(j_{1,6}-1)-1$  \\ \hline
7   & $\frac{\cal P_{7,(1,0)}(7z)-\cal P_{7,(2,0)}(7z)}
            {\cal P_{7,(1,0)}(7z)-\cal P_{7,(4,0)}(7z)}$ &
      $-1/(j_{1,7}-1)-3$  \\ \hline
8   & $\theta_3(2z)/\theta_3(4z)$ & $2/(j_{1,8}-1)-1$ \\ \hline
9   & $\frac{\cal P_{9,(1,0)}(9z)-\cal P_{9,(2,0)}(9z)}
            {\cal P_{9,(1,0)}(9z)-\cal P_{9,(4,0)}(9z)}$ &
      $-1/(j_{1,9}-1)-2$  \\ \hline
10   & $\frac{\cal P_{10,(1,0)}(10z)-\cal P_{10,(2,0)}(10z)}
             {\cal P_{10,(1,0)}(10z)-\cal P_{10,(4,0)}(10z)}$ &
      $-1/(j_{1,10}-1)-2$  \\ \hline
12   & $\theta_3(2z)/\theta_3(6z)$ & $2/(j_{1,12}-1)$ \\ \hline
\end{tabular}
\end{center}
 \vskip0.5cm
 When $N=1,2,3,4,6$, ${\cal N}(j_{1,N})$ becomes a Thompson
 series $T_g$ with $\Gamma_g=\Gamma_0(N)$. Hence, if $N=4$,
 ${\cal N}(j_{1,4})$ has periodically vanishing property
 by Theorem \ref{VVV}-(i).
 Otherwise, the Fourier coefficients of ${\cal N}(j_{1,N})$ do not
 vanish (see \cite{Mckay-Strauss}, Table 1).
 Therefore, we consider only the following cases $N$ for which
 $\overline{\Gamma}_1(N)\neq \overline{\Gamma}_0(N)$.
 \vskip0.5cm
 \begin{center} {\bf Table 4}. Fourier coefficients $H_m$
                       of ${\cal N}(j_{1,N})$ for $1\le m \le 60$
                       \vskip0.2cm
\begin{tabular} {|c||c|c|c|c|c|c|} \hline   & ${\cal N}(j_{1,5})$ & ${\cal N}(j_{1,7})$
                                            & ${\cal N}(j_{1,8})$ & ${\cal N}(j_{1,9})$
                                            & ${\cal N}(j_{1,10})$ & ${\cal N}(j_{1,12})$\\
 \hline\hline
 $H_{1}$ &     10  &       4 &       3 &        2 &        2 &        1 \\ \hline
 $H_{2}$ &      5  &       3 &       2 &        2 &        1 &        1 \\ \hline
 $H_{3}$ &    -15  &       0 &       1 &        1 &        1 &        1 \\ \hline
 $H_{4}$ &    -24  &      -5 &      -2 &       -1 &        0 &        0 \\ \hline
 $H_{5}$ &     15  &      -7 &      -4 &       -2 &       -1 &        0 \\ \hline

 $H_{6}$ &     70  &      -2 &      -4 &       -3 &       -2 &       -1 \\ \hline
 $H_{7}$ &     30  &       8 &       0 &       -2 &       -2 &       -1 \\ \hline
 $H_{8}$ &   -125  &      16 &       6 &        1 &       -1 &       -1 \\ \hline
 $H_{9}$ &   -175  &      12 &       9 &        4 &        1 &       -1 \\ \hline
 $H_{10}$&     95  &      -7 &       8 &        6 &        3 &        0 \\ \hline

 \end{tabular}
 \vskip0.1cm
 \begin{tabular} {|c||c|c|c|c|c|c|} \hline
 & ${\cal N}(j_{1,5})$ & ${\cal N}(j_{1,7})$
                                            & ${\cal N}(j_{1,8})$ & ${\cal N}(j_{1,9})$
                                            & ${\cal N}(j_{1,10})$ & ${\cal N}(j_{1,12})$\\
 \hline\hline

 $H_{11}$&    420  &     -29 &      -1 &        5 &        4 &        1 \\ \hline
 $H_{12}$&    180  &     -35 &     -12 &        1 &        4 &        2 \\ \hline
 $H_{13}$&   -615  &     -10 &     -20 &       -5 &        1 &        2 \\ \hline
 $H_{14}$&   -826  &      37 &     -16 &      -11 &       -2 &        2 \\ \hline
 $H_{15}$&    410  &      70 &       1 &      -12 &       -6 &        1 \\ \hline

 $H_{16}$&   1760  &      53 &      22 &       -7 &       -8 &        0 \\ \hline
 $H_{17}$&    705  &     -21 &      38 &        3 &       -7 &       -2 \\ \hline
 $H_{18}$&  -2415  &    -106 &      30 &       15 &       -3 &       -3 \\ \hline
 $H_{19}$&  -3100  &    -126 &       1 &       22 &        4 &       -4 \\ \hline
 $H_{20}$&   1530  &     -38 &     -40 &       19 &       10 &       -4 \\ \hline

 $H_{21}$&   6270  &     119 &     -64 &        5 &       14 &       -2 \\ \hline
 $H_{22}$&   2460  &     226 &     -52 &      -15 &       12 &        0 \\ \hline
 $H_{23}$&  -8090  &     164 &      -2 &      -32 &        6 &        3 \\ \hline
 $H_{24}$& -10174  &     -70 &      68 &      -36 &       -6 &        5 \\ \hline
 $H_{25}$&   4840  &    -326 &     107 &      -22 &      -16 &        7 \\ \hline

 $H_{26}$&  19570  &    -378 &      88 &        8 &      -22 &        6 \\ \hline
 $H_{27}$&   7500  &    -106 &      -2 &       40 &      -20 &        4 \\ \hline
 $H_{28}$& -24360  &     353 &    -112 &       58 &       -8 &        0 \\ \hline
 $H_{29}$& -30024  &     652 &    -180 &       50 &        8 &       -4 \\ \hline
 $H_{30}$&  14130  &     469 &    -144 &       12 &       26 &       -8 \\ \hline

 $H_{31}$&  55970  &    -189 &       3 &      -41 &       34 &      -10 \\ \hline
 $H_{32}$&  21155  &    -885 &     182 &      -84 &       31 &       -9 \\ \hline
 $H_{33}$& -67380  &   -1015 &     292 &      -93 &       12 &       -6 \\ \hline
 $H_{34}$& -81926  &    -290 &     228 &      -54 &      -14 &        0 \\ \hline
 $H_{35}$&  37895  &     910 &       4 &       22 &      -41 &        6 \\ \hline

 $H_{36}$& 148410  &    1664 &    -286 &      103 &      -54 &       12 \\ \hline
 $H_{37}$&  55305  &    1179 &    -452 &      148 &      -47 &       14 \\ \hline
 $H_{38}$&-174500  &    -483 &    -356 &      124 &      -20 &       14 \\ \hline
 $H_{39}$&-209577  &   -2205 &      -4 &       32 &       23 &        8 \\ \hline
 $H_{40}$&  96025  &   -2492 &     440 &      -96 &       61 &        0 \\ \hline

 $H_{41}$& 371620  &    -692 &     686 &     -200 &       84 &      -10 \\ \hline
 $H_{42}$& 137160  &    2212 &     544 &     -219 &       72 &      -18 \\ \hline
 $H_{43}$&-427665  &    3998 &      -5 &     -128 &       31 &      -22 \\ \hline
 $H_{44}$&-508800  &    2809 &    -668 &       46 &      -32 &      -20 \\ \hline
 $H_{45}$& 230670  &   -1120 &   -1044 &      231 &      -90 &      -12 \\ \hline

 $H_{46}$& 885070  &   -5119 &    -816 &      330 &     -122 &        0 \\ \hline
 $H_{47}$& 323605  &   -5754 &       5 &      275 &     -107 &       15 \\ \hline
 $H_{48}$&-1001340 &   -1598 &     996 &       67 &      -44 &       26 \\ \hline
 $H_{49}$&-1181123 &    4992 &    1563 &     -216 &       45 &       33 \\ \hline
 $H_{50}$& 531545  &    8968 &    1210 &     -439 &      133 &       29 \\ \hline

 \end{tabular}
 \vskip0.1cm
 \begin{tabular} {|c||c|c|c|c|c|c|} \hline
 & ${\cal N}(j_{1,5})$ & ${\cal N}(j_{1,7})$
                                            & ${\cal N}(j_{1,8})$ & ${\cal N}(j_{1,9})$
                                            & ${\cal N}(j_{1,10})$ & ${\cal N}(j_{1,12})$\\
 \hline\hline

 $H_{51}$& 2022670 &    6251 &       6 &     -477 &      174 &       19 \\ \hline
 $H_{52}$&  734130 &   -2506 &   -1464 &     -275 &      154 &        0 \\ \hline
 $H_{53}$&-2253515 &  -11285 &   -2276 &      107 &       61 &      -20 \\ \hline
 $H_{54}$&-2639348 &  -12579 &   -1768 &      501 &      -68 &      -37 \\ \hline
 $H_{55}$& 1178880 &   -3455 &      -8 &      708 &     -192 &      -45 \\ \hline

 $H_{56}$& 4456650 &   10812 &    2128 &      590 &     -254 &      -42 \\ \hline
 $H_{57}$& 1606500 &   19278 &    3284 &      146 &     -220 &      -26 \\ \hline
 $H_{58}$&-4901250 &   13362 &    2552 &     -447 &      -90 &        0 \\ \hline
 $H_{59}$&-5703676 &   -5278 &      -9 &     -911 &      100 &       27 \\ \hline
 $H_{60}$& 2532720 &  -23765 &   -3056 &     -987 &      272 &       52 \\ \hline
 \end{tabular}
 \end{center}

\end{document}